\documentclass{amsart}
\usepackage{amsmath}
\usepackage{amsfonts}
\usepackage{graphics}
\usepackage{epsfig}
\usepackage{amssymb}
\usepackage{amscd}

\newtheorem{thm}{Theorem}[section]
\newtheorem{cor}[thm]{Corollary}
\newtheorem{lem}[thm]{Lemma}
\newtheorem{prop}[thm]{Proposition}
\theoremstyle{definition}
\newtheorem{Def}[thm]{Definition}
\newtheorem*{ack}{Acknowledgement}
\theoremstyle{remark}
\newtheorem{rem}[thm]{Remark}

\numberwithin{equation}{section}
\numberwithin{figure}{section}

\def\n{{\mathbf{n}}}
\def\w{{\mathbf{w}}}
\def\e{{\mathbf{e}}}

\def\isom{\cong}
\def\tensor{\otimes}
\def\dsum{\oplus}
\def\Hom{{\text{\rm{Hom}}}}
\def\Spec{{\text{\rm{Spec}}}}

\def\Proj{{\text{\rm{Proj}}}}
\def\Ker{{\text{\rm{Ker}}}}

\def\Pic{{\text{\rm{Pic}}}}
\def\Coker{{\text{\rm{Coker}}}}
\def\ord{{\text{\rm{ord}}}}
\def\exp{{\text{\rm{exp}}}}
\def\Jac{{\text{\rm{Jac}}}}
\def\Prym{{\text{\rm{Prym}}}}
\def\trace{{\text{\rm{trace}}}}

\def\alp{\alpha}		
		
\def\gam{{\hbox{\raise1pt\hbox{$\gamma$}}}}
		\def\Gam{\Gamma}

\def\lam{\lambda}

\def\rchi{{\hbox{\raise1.5pt\hbox{$\chi$}}}}

\def\longrarrow{\longrightarrow}

\def\hookrarrow{\hookrightarrow}

\begin{document}

\title[Prym varieties and Integrable Systems]
{Prym varieties and Integrable Systems$^\dagger$}

\author[Y.~Li]{Yingchen Li}
\address{
Department of Mathematics\\
University of California\\
Davis, CA 95616--8633}
\author[M.~Mulase]{Motohico Mulase$^\ddagger$}
\address{
Department of Mathematics\\
University of California\\
Davis, CA 95616--8633}
\email{mulase@math.ucdavis.edu}

\allowdisplaybreaks
\setcounter{section}{-1}
\begin{abstract}
A new relation between Prym varieties of
arbitrary morphisms of algebraic curves and integrable systems 
is discovered. The action of
 maximal commutative subalgebras of the formal loop algebra of
$GL_n$ defined on certain
infinite-dimensional Grassmannians is studied.
 It is proved that every finite-dimensional
orbit of the action of traceless elements of these commutative Lie algebras
is  isomorphic to the Prym variety associated with a morphism of 
algebraic curves. Conversely, it is shown that every Prym variety can be
realized as a finite-dimensional orbit of the action of traceless diagonal
elements of the formal loop algebra, which defines the multicomponent
KP system.  
\end{abstract}
\thanks{$^\dagger$Published in Communications in Analysis and Geometry
\textbf{5}, (1997) 279--332. Earlier version was circulated
under the title \emph{Category of morphisms of algebraic curves and a
characterization of Prym varieties} as Max-Planck-Institut Preprint
MPI/92-24.}
\thanks{$^\ddagger$Research 
supported in part by NSF Grant DMS 91--03239, 94--04111.}
\maketitle
\tableofcontents

\section{Introduction}
\label{sec: intro}
\medskip

\noindent\textbf{0.1}
{From} a geometric point of view, the \emph{Kadomtsev-Petviashvili
(KP) equations} are best understood as a set of commuting
vector fields, or \emph{flows}, defined on an infinite-dimensional
Grassmannian \cite{S}. The Grassmannian $Gr_1(\mu)$ 
is the set of vector subspaces
$W$ of the field $L = \mathbb{C}((z))$ of formal Laurent series in $z$
such that the projection $W\longrarrow \mathbb{C}((z))/\mathbb{C}[[z]]z$
is a Fredholm map of index $\mu$. 
The commutative algebra $\mathbb{C}[z^{-1}]$
acts on $L$ by multiplication, and hence it induces 
commuting flows on the Grassmannian. This very simple picture
is nothing but the KP system written in the language
of infinite-dimensional geometry. A striking fact is that 
every finite-dimensional orbit (or integral manifold) of these flows
is canonically isomorphic to the Jacobian variety
of an algebraic curve, and conversely,
every Jacobian variety can be realized as a
finite-dimensional orbit of the KP flows \cite{M1}.
This statement is equivalent to the claim 
that the KP equations
characterize the Riemann theta functions associated with Jacobian
varieties \cite{AD}. 

If one generalizes the above Grassmannian to the Grassmannian
$Gr_n(\mu)$ consisting of vector subspaces of $L^{\dsum n}$ with
a Fredholm condition, then the formal loop algebra 
$gl(n,L)$ acts on it.
In particular, the Borel subalgebra (one of
the maximal commutative subalgebras) of Heisenberg algebras acts
on $Gr_n(\mu)$ with the center acting trivially. Let us
call the system of vector fields coming from this action the 
\emph{Heisenberg flows} on $Gr_n(\mu)$. Now one
can ask a question: what are the finite-dimensional orbits of
these Heisenberg flows, and what kind of geometric
objects do they represent?
Actually, this question was asked to one of the authors 
by Professor H.~Morikawa as early as in 1984. In this paper,
we give a complete answer to this question. Indeed, we shall prove
(see Proposition~\ref{5.1. Proposition} and Theorem~\ref{5.8. Theorem} below)
\medskip
\begin{thm}
\label{A}
A finite-dimensional orbit of the Heisenberg flows defined on
the Grassmannian of vector valued functions corresponds to
a covering morphism of algebraic curves, and 
the orbit itself is canonically isomorphic to the Jacobian
variety of the curve upstairs. Moreover, the action of the 
traceless elements of the Borel subalgebra (the traceless
Heisenberg flows)
produces the Prym variety associated with this covering 
morphism as an orbit.
\end{thm}
\medskip
\begin{rem}
The relation between Heisenberg algebras and covering
morphisms of algebraic curves was first discovered 
by Adams and Bergvelt \cite{AB}.
\end{rem}

\medskip
\noindent\textbf{0.2}
Right after the publication of works (\cite{AD}, \cite{M1}, 
\cite{Sh1}) on characterization of
Jacobian varieties by means of integrable systems,
it has become an important 
problem to find a similar theory for
 Prym varieties. We establish in this paper a simple
solution of this problem in terms of the \emph{multi-component
KP system} defined on a certain quotient space
of the Grassmannian of vector valued
functions.

Classically,  Prym varieties associated with degree two
coverings of algebraic curves were used by Schottky
and Jung in their approach to the Schottky problem.
 The modern interests
 in Prym varieties were revived in \cite{Mum1}. 
  Recently,  Prym varieties of higher degree coverings
have been used in the study of the generalized
 theta divisors on the moduli spaces of stable vector bundles
 over an algebraic curve \cite{BNR}, \cite{H}.
 This direction of research, usually called
 ``Hitchin's Abelianization Program,''  owes its motivation
 and methods to finite dimensional integrable systems in the context
 of symplectic geometry.
In the case  of infinite dimensional  integrable systems, it has been
discovered that Prym varieties of ramified double
sheeted coverings of curves appear as  solutions of the BKP system
 \cite{DJKM}.
Independently, a Prym variety of degree two
covering with exactly two
ramification points has been observed
in the deformation theory of two-dimensional
Schr\"odinger operators  \cite{No}, \cite{NV}.
As far as the authors know, the only
Prym varieties so far
considered in the context of integrable systems
are associated with  ramified, double sheeted
coverings of algebraic curves. Consequently, the attempts (\cite{Sh2},
\cite{T}) of characterizing
Prym varieties in terms of integrable systems are all restricted
to these special Prym varieties.

Let us define the \emph{Grassmannian quotient} $Z_n(0)$ 
as the quotient space of $Gr_n(0)$ by the diagonal action
of $\big(1 + \mathbb{C}[[z]]z\big)^{\times n}$. The \emph{traceless
$n$-component KP system} is defined by the action of the traceless
diagonal matrices with entries in $\mathbb{C}[z^{-1}]$ on $Z_n(0)$.
Since this system is a special case of the traceless Heisenberg
flows, every finite-dimensional orbit of this system is a
Prym variety. 
Conversely, an arbitrary Prym variety associated with a degree $n$
covering morphism of algebraic curves can be realized 
as a finite-dimensional orbit. 
Thus we have (see Theorem~\ref{5.14. Theorem} below):
\medskip
\begin{thm}
\label{B}
An algebraic variety is isomorphic to the Prym variety associated with 
a degree $n$ covering of an algebraic curve if and only if it
can be realized as a finite-dimensional 
 orbit of the traceless $n$-component KP
system defined on the Grassmannian quotient $Z_n(0)$.
\end{thm}

\medskip
\noindent\textbf{0.3}
The relation between algebraic geometry and the Grassmannian
comes from the cohomology map of \cite{SW},
which assigns injectively a point of $Gr_1(0)$ to a 
set of geometric data consisting of an algebraic curve and a
line bundle together with some local information. 
This correspondence was enlarged in \cite{M3} to 
include arbitrary vector bundles on curves.

In this paper, we generalize the cohomology functor
of \cite{M3} so that we can deal with arbitrary  morphisms
between algebraic curves. Let $\n = (n_1, \cdots, n_\ell)$ denote
an integral vector consisting of positive integers satisfying that
$n = n_1 + \cdots + n_\ell$.
\medskip
\begin{thm}
\label{C}
For each $\n$, the following two categories are equivalent:
\begin{enumerate}

\item The category $\mathcal{C}(\n)$. An object of this category
 consists of  an arbitrary degree $n$ morphism
 $f:C_\n\longrarrow C_0$ of algebraic curves 
and an arbitrary vector bundle $\mathcal{F}$ on $C_\n$.
 The curve $C_0$ has a smooth marked point $p$ with a local coordinate
$y$ around it. The curve $C_\n$ has $\ell$ $(1\le \ell \le n)$ 
 smooth marked
points $\{p_1, \cdots, p_\ell\} = f^{-1}(p)$ 
 with ramification index $n_j$ at each point $p_j$.
 The curve $C_\n$ is further endowed
with a local coordinate $y_j$ and a local trivialization
of $\mathcal{F}$ around $p_j$.

\item The category $\mathcal{S}(\n)$. An object of this category
 is  a triple $(A_0,A_\n,W)$  consisting
of a point $W\in \bigcup_{\mu\in \mathbb{Z}}Gr_n(\mu)$, a 
``large'' subalgebra
$A_0\subset \mathbb{C}((y))$ for some $y\in \mathbb{C}[[z]]$, and another
``large'' subalgebra
$$
A_\n\subset \bigoplus_{j=1}^\ell \mathbb{C}((y^{1/n_j}))
\isom \bigoplus_{j=1}^\ell \mathbb{C}((y_j))\;.
$$
In a certain matrix representation as subalgebras of the
formal loop algebra
$gl\big(n,\mathbb{C}((y))\big)$ acting on the Grassmannian, 
they satisfy $A_0\subset A_\n$
and $A_\n\cdot W\subset W$.
\end{enumerate}
\end{thm}
\medskip
\noindent
The precise statement of this theorem is given in Section~\ref{sec: coh}, and 
its proof is completed in Section~\ref{sec: inverse}.
One of the reasons of 
introducing a category rather than just a set is because we
need not only a set-theoretical bijection of objects but
also a canonical correspondence of the morphisms
 in the proof of the claim
that every Prym variety can be realized as a finite-dimensional
 orbit of the traceless multi-component
KP system on the Grassmannian quotient.

\medskip
\noindent\textbf{0.4}
 The motivation of extending the framework of the
original Segal-Wilson construction to include arbitrary vector bundles on curves
of \cite{M3} was to  establish 
a complete geometric classification of
 all the  commutative algebras
consisting of ordinary differential operators with 
coefficients in scalar valued functions. If we 
 apply the functor of Theorem~\ref{C} in this direction,
then we obtain (see Proposition~\ref{6.14. Proposition} and Theorem~\ref{6.15. Theorem} below):
\medskip
\begin{thm}
\label{D}
Every object of the category $\mathcal{C}(\n)$ with a smooth curve
$C_\n$ and a line bundle $\mathcal{F}$ on $C_\n$ satisfying the 
cohomology vanishing condition 
$$
H^0(C_\n,\mathcal{F}) = H^1(C_\n,\mathcal{F}) = 0
$$
gives rise to a maximal commutative algebra consisting of 
ordinary differential operators with coefficients in $n\times n$
matrix valued functions.
\end{thm}
\medskip\noindent
Some examples of commuting matrix ordinary 
differential operators have been studied
 before (\cite{G}, \cite{Na}). Grinevich's work
 is  different from ours. In \cite{G} he considers commuting pairs
 of matrix differential operators. For each commuting pair he constructs a
 single affine algebraic curve (possibly reducible) in the affine plane
  and a vector
bundle on each of the irreducible components and conversely, given such
a collection of algebro-geometric data together with some extra local
information he constructs a commuting pair of matrix differential operators.
In our case, the purpose is to classify commutative algebras of
matrix differential operators. This point of view is more intrinsic than 
considering commuting pairs because they are particular choices of generators
 of the algebras. On the algebro-geometric
side, we obtain morphisms of two abstract curves (no embeddings)
 and maps of the corresponding Jacobian
varieties.  Prym varieties come in very naturally in our picture.
 Nakayashiki's construction (Appendix of \cite{Na}) is
similar to ours, but that corresponds to 
locally cyclic coverings
of curves, i.e.\ a morphism $f:C\longrarrow C_0$ such that
there is a point $p\in C_0$
where $f^{-1}(p)$ consists of one point.
Since we can use arbitrary coverings of curves, 
we obtain in this paper 
 a far larger class of totally new
examples systematically.
As a key step from algebraic geometry of curves 
and vector bundles to the differential operator algebra
with matrix coefficients, we prove the following (see Theorem~\ref{6.5. Theorem} below):

\medskip
\begin{thm}
\label{E}
The big-cell of the Grassmannian $Gr_n(0)$ is 
canonically identified with
the group of monic invertible pseudodifferential operators with 
matrix coefficients.
\end{thm}
\medskip\noindent
Only the case of $n=1$ of this statement was known before. With
this identification, we can translate the flows
on the Grassmannian associated with an arbitrary
commutative subalgebra of the loop algebras
into an integrable system of nonlinear
partial differential equations. The unique solvability of these
systems can be shown  by using the generalized Birkhoff
decomposition of \cite{M2}.

\medskip
\noindent\textbf{0.5}
 This paper is organized as follows. In Section~\ref{sec: 1}, we review
 some  standard facts about Prym varieties. The Heisenberg flows are
introduced in Section~\ref{sec: Heis}. Since we do not deal with any central 
extensions in this paper, we shall not use the Heisenberg algebras 
in the main text. All we need are the maximal commutative
subalgebras of the formal loop algebras. Accordingly, the action
of the Borel subalgebras will be replaced by the action of
the full maximal commutative algebras defined on certain
quotient spaces of the Grassmannian. This turns out to be
more natural because of the coordinate-free nature of the flows
on the quotient spaces. The two categories we work with are defined
in Section~\ref{sec: coh}, where a generalization of the cohomology
functor is given. 
In Section~\ref{sec: inverse}, we give the construction of the geometric
data out of the algebraic data consisting of commutative
algebras and a point of the Grassmannian.
The finite-dimensional orbits of the Heisenberg flows
are studied in Section~\ref{sec: Prym}, in which the characterization
theorem of Prym varieties is proved.
Section~\ref{sec: ODE} is devoted to explaining the relation of the entire
theory with the ordinary differential operators with 
matrix coefficients. 

The results we obtain in Sections~\ref{sec: coh}, \ref{sec: inverse},
 and \ref{sec: ODE} (except for Theorem~\ref{6.15. Theorem}, where we need
zero characteristic) hold for an arbitrary field $k$.
In Sections~\ref{sec: 1} and \ref{sec: Prym} (except for Proposition~\ref{5.1. Proposition}, which is true
for any field), we work with the field $\mathbb{C}$ of complex
numbers.
\medskip

\begin{ack}
The authors wish to express their gratitude 
to the Max-Planck-Institut f\"ur Mathematik
for generous financial support and hospitality, without which the entire
project would never have  taken place. They also thank 
 M.~Bergvelt
for sending them the paper \cite{AB} prior to its publication. The earlier version of the current article has been circulated as a Max-Planck-Institut preprint.
\end{ack}

\bigskip

\section{Covering morphisms of curves and  Prym varieties}
\label{sec: 1}
\medskip

We begin with defining Prym varieties in the most general
setting, and then introduce \emph{locally cyclic coverings} of curves,
which play an important role in defining the category of
arbitrary covering morphisms of algebraic curves in Section~\ref{sec: coh}.

\begin{Def} \label{1.1. Definition} Let
$f : C\longrightarrow {C_0}$ be a covering morphism of degree $n$
between smooth algebraic curves
$C$ and ${C_0}$, and let
$N_f: \Jac(C)\longrarrow \Jac({C_0})$ be the norm homomorphism between
the Jacobian varieties,
which assigns to an element $\sum_q n_q\cdot q\in \Jac(C)$
its image $\sum_q n_q\cdot f(q)\in \Jac({C_0})$. This
 is a surjective homomorphism, and hence the kernel 
$\Ker(N_f)$
   is an abelian subscheme of $\Jac(C)$ of dimension
   $g(C)-g({C_0})$, where $g(C)$ denotes the genus of the curve $C$.
  We call this kernel the \emph{Prym variety}
 associated  with the morphism  $f$,  and denote it by $\Prym(f)$.
\end{Def}
\medskip

\begin{rem} \label{{1.2.} Remark} Usually the Prym variety
 of a covering morphism $f$ is defined to be the
 connected component of the kernel of the
 norm homomorphism containing 0. 
Since any two connected components of $\Ker(N_f)$
 are translations of each other in $\Jac(C)$, there is no harm to
 call the whole kernel the Prym variety.
   If  the pull-back homomorphism
 $f^{*}: \Jac({C_0})\longrarrow \Jac(C)$ is injective, then  
  the norm homomorphism can be identified with the transpose
 of $f^{*}$, and hence its kernel is connected. So in this situation,
 our definition coincides with the usual one. We will give 
a class of coverings 
  where the pull-back homomorphisms  are injective (see Proposition~\ref{1.7. Proposition}). 
\end{rem}

\medskip
\begin{rem} \label{{1.3.} Remark} 
  Let ${R}\subset C$ be the ramification divisor  of the morphism
 $f$ of Definition~\ref{1.1. Definition} and $\mathcal{O}_C({R})$ the locally free sheaf
associated with ${R}$. Then it can be shown that for any
 line bundle $\mathcal{L}$ on $C$, we have  $N_f(\mathcal{L})=\det(f_*\mathcal{L})\otimes
    \det\big(f_*\mathcal{O}_C({R})\big)$.
   Thus up to a translation, the norm homomorphism can be
 identified with the map assigning the determinant of the direct
image to the line bundle on $C$. Therefore, one can talk about the
 Prym varieties in $\Pic^d(C)$ for an arbitrary $d$,
 not just in $\Jac(C)=\Pic^{0}(C)$.

When the curves $C$ and $C_0$ are singular, we replace the 
Jacobian variety $\Jac(C)$ by the generalized Jacobian,
which is the connected component of $H^1(C, \mathcal{O}^*_{C})$
containing the structure sheaf. By taking the determinant 
of the direct image sheaf, we can define a map of the
generalized Jacobian of $C$ into $H^1(C_0,\mathcal{O}^*_{C_0})$.
The fiber of this map is called the \emph{generalized} Prym variety
associated with the morphism $f$.
\end{rem}
\medskip
\begin{rem} \label{{1.4.} Remark} 
    According to our definition (Definition~\ref{1.1. Definition}), the  Jacobian variety
 of an arbitrary algebraic curve $C$ can
 be viewed as a Prym variety. Indeed, for a nontrivial
 morphism of $C$ onto $\mathbb{P}^1$,  the induced
 norm homomorphism is the zero-map.   Thus the class of Prym varieties
 contains  Jacobians as a subclass. Of course there are 
 infinitely many ways to realize $\Jac(C)$ as a Prym variety
 in this manner. 

\end{rem}
\medskip

Let us consider the polarizations  of  Prym varieties. 
Let $\Theta_{C}$ and $\Theta_{{C_0}}$ be the Riemann
 theta divisors on $\Jac(C)$ and $\Jac({C_0})$, respectively. Then
  the restriction of $\Theta_{C}$ to $\Prym(f)$ gives an
 ample divisor $H$  on $\Prym(f)$. This is usually not a 
 principal polarization if $g(C_0)\not= 0$. 
There is a natural homomorphism 
$\psi:\Jac({C_0})\times \Prym(f)\longrarrow \Jac(C)$
  which assigns $f^*\mathcal{L}\otimes \mathcal{M}$
to $(\mathcal{L},\mathcal{M})\in \Jac({C_0})\times \Prym(f)$. This is
 an isogeny, and the pull-back of $\Theta_{C}$ under this homomorphism
 is given by 
$$
\psi^*\mathcal{O}_{\Jac(C)}(\Theta_{C}) 
\isom \mathcal{O}_{\Jac({C_0})}(n\Theta_{C_0})\otimes \mathcal{O}_{\Prym(f)}(H)\;.
$$

In Section~\ref{sec: coh}, we define a category of covering morphisms
of algebraic curves. As a \emph{morphism} between the 
covering morphisms, we use the following special coverings:
\medskip
\begin{Def} \label{1.5. Definition} A degree $r$  morphism 
$\alpha : C\longrarrow {C_0}$ of 
   algebraic curves is said to
be a \emph{locally cyclic covering} if there
is a point $p\in {C_0}$ such that $\alpha^{*}(p)=r\cdot q$ for some 
$q\in C$.
\end{Def}
\medskip
\begin{prop} \label{1.6. Proposition}  Every smooth projective curve $C$
 has  infinitely many smooth
 locally cyclic coverings of an arbitrary degree.
\end{prop}
\medskip
\begin{proof}  We use the theory of spectral curves
 to prove this statement. For a detailed account of
 spectral curves, we refer to \cite{BNR} and \cite{H}.

Let us take a line bundle $\mathcal{L}$ over $C$ of sufficiently
 large degree. For such $\mathcal{L}$ we can choose sections
  $s_{i}\in H^{0}(C, \mathcal{L}^i)$, $i$ $=$ 1, 2, $\cdots$, $r$, 
satisfying the following conditions:
\begin{enumerate}
\item  All  $s_{i}$'s have a common zero point, say $p\in C$,
 i.e., $s_i\in H^0(C, \mathcal{L}^i(-p))$, ${\ } {\ } {\ }$
 $i=1,2,\cdots , r$;
\item  $s_{r}\notin  H^{0}(C, \mathcal{L}^{r}(-2p))$.
\end{enumerate}
    Now consider the sheaf $\mathcal{R}$ of symmetric  $\mathcal{O}_{C}$-algebras
generated by $\mathcal{L}^{-1}$. As an $\mathcal{O}_{C}$-module
 this algebra can be written as 
$$
\mathcal{R} = \bigoplus_{i=0}^{\infty }\mathcal{L}^{-i}\;.
$$
In order to construct a locally cyclic covering of $C$,
we take the ideal $\mathcal{I}_s$
 of the algebra $\mathcal{R}$ generated by
 the image of the sum of the homomorphisms 
$s_{i}: \mathcal{L}^{-r}\longrarrow \mathcal{L}^{-r+i}$. 
We define $C_{s}=\Spec(\mathcal{R}/\mathcal{I}_s)$, where
$s = (s_1, s_2, \cdots, s_r)$. Then $C_{s}$ is
a  spectral curve, and  
the natural projection $\pi:C_s\longrarrow C$ gives a
degree $r$ covering of $C$.  For sufficiently general
 sections $s_{i}$ with properties (1) and (2), we may also assume 
the following (see \cite{BNR}):
\begin{enumerate}
\item[(3)] The spectral curve  $C_{s}$ is integral, 
i.e.~reduced and irreducible.
\end{enumerate}
We claim here that  $C_{s}$ is smooth
 in a neighborhood of the inverse image of $p$. 
Indeed,  let us take a local
 parameter $y$ of $C$ around $p$ and a local coordinate $x$ in the
fiber direction of the total space of the line bundle $\mathcal{L}$. 
Then the local Jacobian criterion 
for smoothness
 in a neighborhood of $\pi^{-1}(p)$ states that 
the following system 
$$
\begin{cases}
x^{r} +s_{1}(y)x^{r-1}+\cdots +s_{r}(y)=0\\
rx^{r-1}+s_{1}(y)(r-1)x^{r-2} +\cdots +s_{r-1}(y)=0\\
s_{1}(y)'x^{r-1}+ s_{2}(y)'x^{r-2}+\cdots +s_{r}(y)'=0
\end{cases}
$$
of equations in $(x, y)$
 has no solutions.  But  this is  clearly the case in our situation
because of the conditions (1), (2) and (3).
Thus we have verified the claim.
It is also clear that $\pi^*(p) = r\cdot q$, where $q$ is
the point of $C_s$ defined by $x^r = 0$ and $y = 0$.
 Then by taking the
 normalization 	of $C_{s}$ we obtain a smooth
 locally cyclic covering of $C$. This completes the proof.
\end{proof}
\medskip
\begin{prop} \label{1.7. Proposition}
Let $\alpha:C\longrarrow {C_0}$ be a locally cyclic covering of
degree $r$. Then the induced
homomorphism $\alpha^*:\Jac({C_0})\longrarrow \Jac(C)$
of Jacobians is injective. In particular, the Prym variety 
$\Prym(\alp)$
associated with the morphism $\alp$ is
connected.
\end{prop}
\medskip
\begin{proof}
  Let us suppose in contrary that
 $\mathcal{L}\not\cong {\mathcal{O}}_{{C_0}}$ and  
$\alpha^{*}\mathcal{L}\cong {\mathcal{O}}_{C}$ 
for some 
 $\mathcal{L}\in \Jac({C_0})$. Then by the projection formula we have
 $ \mathcal{L}\otimes \alpha_{*}{\mathcal{O}}_{C}\cong \alpha_{*}\mathcal{O}_{C}$.
 Taking determinants on both sides we see
   that $\mathcal{L}$ is an $r$-torsion point in $\Jac({C_0})$,
i.e.~$\mathcal{L}^r \cong \mathcal{O}_{C_0}$.  Let $m$ be the smallest positive
 integer satisfying that ${\mathcal{L}}^{m} \cong {\mathcal{O}}_{{C_0}}$.
 Let us consider the spectral curve 
$$
C' = \Spec\big( \bigoplus _{i=0}^{\infty} \mathcal{L}^{-i}/\mathcal{I}_s \big)
$$ 
defined by the line bundle
 $\mathcal{L}$ and its sections 
$$
s = (s_1, s_2, \cdots, s_{m-1}, s_m) = (0, 0, \cdots, 0, 1)
\in \bigoplus_{i = 1}^m H^{0}({C_0},\mathcal{L}^i)\;.
$$
It is easy to verify that $C'$ is an unramified
 covering of ${C_0}$ of degree $m$. 
Now we claim  that
 the morphism $\alpha : C\longrarrow {C_0}$ factors through $C'$, 
but this leads to a contradiction to our assumption that $\alpha$ is a 
locally cyclic covering.

The construction of such a morphism $f:C\longrarrow C'$
over ${C_0}$ amounts to defining an
 ${\mathcal{O}}_{{C_0}}$-algebra homomorphism  
\begin{equation}
f^{\sharp} :\bigoplus _{i=0}^{\infty}\mathcal{L}^{-i} /\mathcal{I}_s
\longrightarrow 
 \alpha_{*}{\mathcal{O}}_{C}\;.\label{1.8}
\end{equation}
In order to give (\ref{1.8}), it is sufficient to 
define an ${\mathcal{O}}_{{C_0}}$-module homomorphism
  $\phi :\mathcal{L}^{-1}\longrarrow \alpha_{*}{\mathcal{O}}_{C}$ such
  that $\phi ^{\otimes m}: \mathcal{L}^{-m}\cong \mathcal{O}_{C_0}\longrightarrow
 \alpha _{*}\mathcal{O}_{C}$ is the inclusion map induced by $\alpha$.
Since we have
$$
H^0(C,\mathcal{O}_C) \isom 
H^0(C_0,\alpha_*\mathcal{O}_C)
\isom H^{0}(C_0, \mathcal{L}\otimes \alpha _{*}\mathcal{O}_{C})
\isom H^{0}(C_0, \mathcal{L}^{m}\otimes \alpha _{*}\mathcal{O}_{C})\;,
$$
 the existence of the desired $\phi$ is obvious.
 This completes the proof.
\end{proof}

\bigskip

\section{The Heisenberg flows on the Grassmannian of 
vector valued functions}
\label{sec: Heis}
\medskip

In this section, we define the Grassmannians of vector valued
functions and introduce various vector fields (or flows) 
on them. Let $k$ be an arbitrary field, $k[[z]]$ the ring of
formal power series in one variable $z$ defined over $k$, 
and $L = k((z))$ the field of fractions of $k[[z]]$. An element of
$L$ is a formal Laurent series in $z$ with a pole of finite order. 
We call $y = y(z)\in L$ an element of \emph{order} $m$ if 
$y\in k[[z]]z^{-m}\setminus k[[z]]z^{-m+1}$.
Consider the infinite-dimensional
vector space $V = L^{\dsum n}$ over $k$.
It has a natural filtration by the (pole) order
$$
\cdots \subset F^{(m-1)}(V) \subset F^{(m)}(V) \subset F^{(m+1)}(V)\subset 
\cdots\;,
$$
where we define
\begin{equation}
F^{(m)}(V) = \left\{\left. \sum_{j = 0}^\infty a_j z^{-m+j} \;
\right|\; a_j\in  k^{\dsum n}\right\}\;.\label{2.1}
\end{equation}
In particular, we have $F^{(m)}(V)\big/F^{(m-1)}(V)
\isom k^{\dsum n}$ for all $m\in \mathbb{Z}$. The filtration satisfies
$$
\bigcup_{m=-\infty}^\infty F^{(m)}(V) = V\quad{\text{   and   }}\quad
\bigcap_{m=-\infty}^\infty F^{(m)}(V) = \{0\}\;,
$$
and hence it determines a topology in $V$. In Section~\ref{sec: inverse}, 
we will introduce other filtrations of $V$ in order to define
algebraic curves and vector bundles on them. 
The current filtration (\ref{2.1}) is used only
for the purpose of defining the Grassmannian as a pro-algebraic
variety (see for example \cite{KSU}).
\medskip
\begin{Def} \label{2.2. Definition} For every integer $\mu$,
the following set of vector subspaces $W$ of $V$
is called the index $\mu$ Grassmannian 
of vector valued functions
of size $n$:
$$
Gr_n(\mu) = \{ W\subset V\;|\; \gamma_W 
{\text{  is Fredholm of 
index  }} \mu\}\;,
$$
where $\gamma_W:W\longrarrow V\big/F^{(-1)}(V)$ 
is the natural projection.
\end{Def}
\medskip\noindent
Let $N_W = \big\{\ord_z(v) \; \big|\; v\in W\big\}$. Then
the Fredholm condition implies that $N_W$ is bounded from 
below and contains all sufficiently large positive
integers. But of course, this condition of $N_W$ does not imply
the Fredholm property of $\gam_W$ when $n>1$.

\medskip

\begin{rem} \label{{2.3.} Remark}
We have used $F^{(-1)}(V)$ in the above definition as a reference
open set for the Fredholm condition. This is because it becomes
the natural choice in Section~\ref{sec: ODE} when we deal with the differential
operator action on the Grassmannian. {From} purely algebro-geometric
point of view, $F^{(0)}(V)$ can also be used (see Remark~\ref{{4.6.} Remark}).
\end{rem}
\medskip\noindent
The \emph{big-cell} $Gr_n^+(0)$ of the Grassmannian 
of vector valued functions of size $n$ is the set of
vector subspaces $W\subset V$ such that $\gamma_W$ is an isomorphism.
For every point $W\in Gr_n(\mu)$, the tangent space
at $W$ is naturally identified with the space of continuous
homomorphism of $W$ into $V/W$:
$$
T_WGr_n(\mu) = \Hom_{\text{cont}}(W,V/W)\;.
$$
Let us define various vector fields on the Grassmannians.
Since the formal loop algebra $gl(n,L)$ acts on $V$,
every element $\xi\in gl(n,L)$ defines a homomorphism
\begin{equation}
W\longrarrow V\overset{\xi}{\longrarrow} V\longrarrow V/W\;,\label{2.4}
\end{equation}
which we shall denote by ${\Psi}_W(\xi)$. Thus the association
$$
Gr_n(\mu)\owns W\longmapsto {\Psi}_W(\xi)\in \Hom_{\text{cont}}(W,V/W) 
= T_WGr_n(\mu)
$$
determines a vector field ${\Psi}(\xi)$ on the Grassmannian.
For a subset $\Xi\subset gl(n,L)$, we use the notations
${\Psi}_W(\Xi) = \big\{{\Psi}_W(\xi)\;\big|\;\xi\in\Xi\big\}$ and 
${\Psi}(\Xi) = \big\{{\Psi}(\xi)\;\big|\;\xi\in\Xi\big\}$.
\medskip
\begin{Def} \label{2.5. Definition}
A smooth subvariety $X$ of $Gr_n(\mu)$ is said to be an {\it
orbit} (or the \emph{integral manifold}) of the flows of
${\Psi}(\Xi)$ if the tangent space $T_WX$ of $X$ at $W$ is equal to
${\Psi}_W(\Xi)$ as a subspace of the whole
tangent space $T_WGr_n(\mu)$ for every point $W\in X$.
\end{Def}
\medskip
\begin{rem} \label{{2.6.} Remark}
There is a far larger algebra than the loop algebra,
the algebra $gl(n,E)$ of \emph{pseudodifferential operators} with
matrix coefficients, acting on $V$. We will come back to this point
in Section~\ref{sec: ODE}.
\end{rem}
\medskip

Let us choose a monic element
\begin{equation}
y = z^r + \sum_{m=1}^\infty c_mz^{r+m} \in L \label{2.7}
\end{equation}
of order $-r$ and consider the following $n\times n$ matrix
\begin{equation}
h_n(y) = 
\begin{pmatrix}
0 &   &      &   & 0 & y\\
1 & 0 &      &   &   & 0\\
  & 1 &\ddots\\
  &   &\ddots& 0\\
  &   &      & 1 &  0\\
  &   &      &   &  1 & 0
\end{pmatrix}\label{2.8}
\end{equation}
satisfying that $h_n(y)^n = y\cdot I_n$, where
$I_n$ is the identity matrix of size $n$. We denote by
$H_{(n)}(y)$ the algebra generated by $h_n(y)$ over $k((y))$,
which is a maximal commutative subalgebra of the formal loop algebra
$gl\big(n,k((y))\big)$. Obviously, we have a natural $k((y))$-algebra
isomorphism 
$$
H_{(n)}(y) \isom k((y))[x]/(x^n - y)\isom k((y^{1/n}))\;,
$$
where $x$ is an indeterminate.
\medskip
\begin{Def} \label{2.9. Definition}
For every integral vector $\n = (n_1, n_2, \cdots, n_\ell)$
of positive integers $n_j$ such that
$n=n_1 + n_2 + \cdots + n_\ell$ and a monic element $y\in L$ of
order $-r$, we define a maximal commutative $k((y))$-subalgebra
of $gl\big(n,k((y)\big)$ by
$$
H_\n(y) = \bigoplus_{j=1}^\ell H_{(n_j)}(y) \isom 
\bigoplus_{j=1}^\ell k((y^{1/n_j}))\;,
$$
where each $H_{(n_j)}(y)$ is embedded 
by the disjoint principal diagonal blocks:
$$
\begin{pmatrix}
H_{(n_1)}(y)\\
&H_{(n_2)}(y)\\
&&\ddots\\
&&&H_{(n_\ell)}(y)
\end{pmatrix}\;.
$$
The algebra $H_\n(y)$ is called
the \emph{maximal commutative algebra of type}
$\n$ associated with the variable $y$. 
\end{Def}
\medskip\noindent
As a module over the field $k((y))$, the algebra $H_\n(y)$ has  
dimension $n$.
\medskip
\begin{rem} \label{{2.10.} Remark}
The lifting of the algebra $H_\n(y)$ to the
central extension of the formal loop algebra $gl\big(n,k((y))\big)$
is the Heisenberg algebra associated with the conjugacy class
of the Weyl group of $gl(n,k)$ determined by the integral vector $\n$
(\cite{FLM}, \cite{Ka}, \cite{PS}).
The word \emph{Heisenberg} in the following definition has 
its origin in this context.
\end{rem}
\medskip
\begin{Def} \label{2.11. Definition}
The set of commutative vector fields ${\Psi}(H_\n(y))$
defined on $Gr_n(\mu)$ is called 
the \emph{Heisenberg flows}
of type $\n = (n_1, n_2, \cdots, n_\ell)$ and rank $r$
associated with the algebra $H_\n(y)$ and the
coordinate $y$ of $(\ref{2.7})$. Let $H_\n(y)_0$ denote the 
subalgebra of $H_\n(y)$ consisting of the 
traceless elements. The system of vector fields 
${\Psi}\big(H_\n(y)_0\big)$ is called the \emph{traceless Heisenberg
flows}.
The set of commuting vector fields ${\Psi}\big(k((y))\big)$ on 
$Gr_n(\mu)$ is called the \emph{$r$-reduced} KP system 
(or the \emph{$r$-reduction} of the KP system) associated with 
the coordinate $y$.
The usual KP system is defined to be the $1$-reduced KP system with the 
choice of $y = z$. The Heisenberg flows associated with
$H_{(1,\cdots,1)}(z)$ of type $(1, \cdots, 1)$ is called the
\emph{$n$-component} KP system. 
\end{Def}
\medskip
\begin{rem} \label{{2.12.} Remark}
As we shall see in Section~\ref{sec: inverse}, the $H_\n(y)$-action on $V$ is 
equivalent to the component-wise multiplication of (\ref{4.1}) to
(\ref{4.4}). {From} this point of view, the Heisenberg flows of
type $\n$ and rank $r$ are contained in
 the $\ell$-component KP system. What is important in 
our presentation as the Heisenberg flows is the new
algebro-geometric interpretation of the orbits of 
these systems defined on the (quotient) Grassmannian
which can be seen only through the right choice of the
coordinates.
\end{rem}
\medskip
\begin{rem} \label{{2.13.} Remark}
The traceless Heisenberg flows of type $\n = (2)$ and rank one
are known to be equivalent to the BKP system. As we shall see
later in this paper, these flows produce the Prym variety
associated with a double sheeted covering of algebraic curves
with at least one ramification point. This explains why the
BKP system is related only with these very special Prym varieties.
\end{rem}
\medskip
\noindent
The flows defined above are too large from the geometric point of view.
The action of the negative order elements of $gl(n,L)$ should be
considered trivial in order to give a direct connection between
the orbits of these flows and the Jacobian varieties.
Thus it is more convenient to define these flows on 
certain quotient spaces.
So let
\begin{equation}
H_\n(y)^- = H_\n(y)\cap gl\big(n,k[[y]]y\big)
\label{2.14}
\end{equation}
and define an abelian group
\begin{equation}
\Gam_\n(y) =  \exp\big(H_\n(y)^-\big) = 
I_n + H_{\n}(y)^-\;.\label{2.15}
\end{equation}
This group is isomorphic to an affine space, and 
acts on the Grassmannian without fixed points. This can be 
 verified as follows. Suppose we have
  $g\cdot W= W$ for some $g =I_n + h\in \Gamma_\n(y)$
 and $W\in Gr_\n(\mu)$. Then $h\cdot W\subset W$.
  Since  $h$ is a nonnilpotent element of negative order,
 by iterating the action
 of $h$ on $W$, we get a contradiction to  the Fredholm
 condition of $\gam_W$.

\medskip
\begin{Def} \label{2.16. Definition}
The \emph{Grassmannian quotient} of type $\n$, index $\mu$
and rank $r$ associated with the algebra $H_\n(y)$
is the quotient space 
$$
Z_\n(\mu,y) = Gr_n(\mu)\big/\Gam_\n(y)\;.
$$
We denote by $Q_{\n,y}: Gr_n(\mu)\longrarrow Z_\n(\mu, y)$
the canonical projection.
\end{Def}
\medskip\noindent
Since $\Gam_\n(y)$ is an affine space acting on the 
Grassmannian without fixed points, the affine principal
fiber bundle $Q_{\n,y}: Gr_n(\mu)\longrarrow Z_\n(\mu, y)$
is trivial. If the Grassmannian is modeled on a complex
Hilbert space, then one can introduce a K\"ahler structure on it, 
which gives rise to a canonical connection on the principal
bundle $Q_{\n,y}$. In that case, there is a standard way of
defining vector fields on the Grassmannian quotient by
using the connection. In our case, however,
since the Grassmannian $Gr_n(\mu)$ is modeled over
$k((z))$, we cannot use these technique of infinite-dimensional
complex geometry. Because of this reason, instead of
defining vector fields on the Grassmannian quotient, 
we give directly a definition of orbits on $Z_\n(\mu, y)$ in the
following manner.

\medskip
\begin{Def} \label{2.17. Definition}
A subvariety $\overline{X}$ of the quotient
Grassmannian $Z_\n(\mu,y)$ is said to be an \emph{orbit} of the
Heisenberg flows associated with $H_\n(y)$ if 
the pull-back $Q_{\n,y}^{-1}(\overline{X})$ is an orbit of the
Heisenberg flows on the Grassmannian $Gr_n(\mu)$.
\end{Def}
\medskip\noindent
Here, we note that 
because of the commutativity of the algebra $H_\n(y)$ and the
group $\Gam_\n(y)$, the Heisenberg flows on the
Grassmannian ``descend'' to the Grassmannian quotient. 
Thus for the flows generated by subalgebras of $H_\n(y)$, 
we can safely talk about the \emph{induced flows} on the
Grassmannian quotient.
\medskip

\begin{Def} \label{2.18. Definition}
An orbit $X$ of the vector fields $\Psi(\Xi)$ on the Grassmannian
$Gr_n(\mu)$ is said to be  of 
\emph{finite type} if $\overline{X} = 
Q_{\n,y}(X)$ is a finite-dimensional
subvariety of the Grassmannian quotient $Z_\n(\mu,y)$.
\end{Def}
\medskip

In Section~\ref{sec: Prym}, we study algebraic geometry of
finite type orbits of the Heisenberg
flows and establish a characterization of Prym varieties
in terms of these flows.
The actual system of nonlinear partial differential equations 
corresponding to these
vector fields are derived in Section~\ref{sec: ODE}, where the unique
solvability of the initial value problem of these
nonlinear equations is shown by using a theorem of \cite{M2}.

\bigskip
\section{The cohomology functor for covering morphisms of 
algebraic curves}
\label{sec: coh}
\medskip

Krichever \cite{Kr} gave a construction
of an exact solution of the entire KP system out of
a set of algebro-geometric data consisting of curves and line
bundles on them. This construction was formulated as a map
of the set of these geometric data into the Grassmannian by
Segal and Wilson \cite{SW}. Its generalization to the
geometric data containing arbitrary vector bundles on curves
was discovered in \cite{M3}. In order to deal with arbitrary
covering morphisms of algebraic curves, we have to enlarge
the framework of the cohomology functor of \cite{M3}. 

\medskip
\begin{Def} \label{3.1. Definition}
A set of \emph{geometric data} of a covering morphism
of algebraic curves
of type $\n$, index $\mu$
and rank $r$ is the collection 
$$
\left\langle f:\big( C_\n, \Delta, \Pi, \mathcal{F}, \Phi\big)
\longrarrow \big({C_0}, p, \pi, f_*\mathcal{F}, \phi\big)\right\rangle
$$
of the following objects:
\begin{enumerate}

\item $\n = (n_1, n_2, \cdots, n_\ell)$ is an integral vector
 of  positive
integers $n_j$ such that $n = n_1 + n_2 + \cdots + n_\ell$.

\item $C_\n$ is a reduced algebraic curve defined over $k$, and
$\Delta = \{p_1,p_2,\cdots,p_\ell\}$ is a set of 
$\ell$ smooth rational points of $C_\n$.

\item $\Pi = (\pi_1, \cdots, \pi_\ell)$ consists of a cyclic
covering morphism $\pi_j:U_{oj}\longrarrow U_j$ of
degree $r$ which maps the formal completion $U_{oj}$ of the affine line
$\mathbb{A}^1_k$ along the origin onto the formal completion $U_j$ of the
curve ${C_\n}$ along $p_j$.

\item $\mathcal{F}$ is a torsion-free sheaf of rank $r$ defined over $C_\n$
satisfying that
$$
\mu = \dim_k H^0(C_\n,\mathcal{F}) - \dim_k H^1(C_\n,\mathcal{F})\;.
$$

\item $\Phi = (\phi_1,\cdots, \phi_\ell)$ consists of an
$\mathcal{O}_{U_j}$-module isomorphism 
$$
\phi_j:\mathcal{F}_{U_j}\overset{\sim}{\longrarrow}
\pi_{j*}\big(\mathcal{O}_{U_{oj}}(-1)\big)\;,
$$
where $\mathcal{F}_{U_j}$ is the formal completion of $\mathcal{F}$ along $p_j$.
We identify $\phi_j$ and $c_j\cdot \phi_j$ for every nonzero
constant $c_j\in k^*$.

\item $C_0$ is an integral curve with a marked smooth rational 
point $p$.

\item $f:C_\n\longrarrow {C_0}$ is a finite 
morphism of degree $n$ of $C_\n$
onto ${C_0}$ such that $f^{-1}(p) = \{p_1,\cdots,p_\ell\}$
with ramification index $n_j$ at each point $p_j$.

\item $\pi: U_o\longrarrow U_p$ 
is a cyclic covering morphism of
degree $r$  which maps the formal completion 
$U_{o}$ of the affine line
$\mathbb{A}^1_k$ at the origin onto the formal completion $U_p$ of the
curve ${C_0}$ along $p$.

\item $\pi_j:U_{oj}\longrarrow U_j$ and the formal completion
$f_j:U_j\longrarrow U_p$ of the morphism $f$ at $p_j$ satisfy the
commutativity of the diagram
$$
\begin{CD}
U_{oj} @>{\pi_j}>> U_j\\
@V{\psi_j}VV @VV{f_j}V\\
U_o @>>{\pi}> U_p,
\end{CD}
$$
where $\psi_j:U_{oj}\longrarrow U_o$ is a cyclic covering of degree
$n_j$.

\item $\phi:(f_*\mathcal{F})_{U_p}\overset{\sim}{\longrarrow}
\pi_*\bigg(\bigoplus_{j=1}^\ell \psi_{j*}\big(\mathcal{O}_{U_{oj}}(-1)
\big)\bigg)$ is an 
$\big(f_*\mathcal{O}_{C_\n}\big)_{U_p}$-module isomorphism
of the sheaves on the formal scheme $U_p$ which is compatible with
the datum $\Phi$ upstairs.
\end{enumerate}
\end{Def}
\medskip\noindent
Here we note that we have an isomorphism
$\psi_{j*}\big(\mathcal{O}_{U_{oj}}(-1)
\big) \isom \mathcal{O}_{U_o}(-1)^{\dsum n_j}$
as an $\mathcal{O}_{U_o}$-module.

Recall that the original cohomology functor is really a cohomology
functor. In order to see what kind of algebraic data come up
from our geometric data, let us apply the cohomology functor to them.
We choose a coordinate $z$ on the formal scheme 
$U_o$ and fix it
once for all. Then we have $U_o = \Spec\big(k[[z]]\big)$.
Since $\psi_j:U_{oj}\longrarrow U_o$ is a cyclic covering
of degree $n_j$, we can identify $U_{oj} =  
\Spec\big(k[[z^{1/n_j}]]\big)$ so that $\psi_j$ is given by
$z = \big(z^{1/n_j}\big)^{n_j} = z_j^{n_j}$, where $z_j = 
z^{1/n_j}$ is a coordinate of $U_{oj}$.
The morphism $\pi$ determines a coordinate
$$
y = z^r + \sum_{m = 1}^\infty c_m z^{r+m}
$$
on $U_p$.  We also choose a coordinate 
$y_j = y^{1/n_j}$ of $U_j$ in which the morphism $f_j$ 
can be written as 
$y = \big(y^{1/n_j}\big)^{n_j} = y_j^{n_j}$.
Out of the geometric data, we can assign a vector
subspace $W$ of $V$ by
\begin{equation}
\begin{aligned}
W &= \phi\big(H^0({C_0}\setminus\{p\},f_*\mathcal{F})\big)\\
&\subset H^0\bigg(U_p\setminus\{p\},
\pi_*\bigoplus_{j=1}^\ell \psi_{j*}\big(\mathcal{O}_{U_{oj}}(-1)\big) \bigg)\\
&= H^0\bigg(U_o\setminus\{o\},
\bigoplus_{j=1}^\ell \psi_{j*}\big(\mathcal{O}_{U_{oj}}(-1)\big)\bigg)\\
&\isom H^0\big(U_o\setminus\{o\},
\bigoplus_{j=1}^\ell \mathcal{O}_{U_o}(-1)^{\dsum n_j}\big)\\
&\isom H^0\big(U_o\setminus \{o\}, \mathcal{O}_{U_o}(-1)^{\dsum n}\big)
= k((z))^{\dsum n} = V\;.
\end{aligned}\label{3.2}
\end{equation}
Here, we have used the convention of \cite{M3} that 
$$
\begin{aligned}
H^0(C_0\setminus \{p\},\mathcal{O}_{C_0}) &= \varinjlim_m
H^0\big(C_0, \mathcal{O}_{C_0}(m\cdot p)\big)\\
H^0(U_o\setminus \{o\},\mathcal{O}_{U_o}) &= \varinjlim_m
H^0\big(U_o, \mathcal{O}_{U_o}(m)\big) = k((z))\;,
\end{aligned}
$$
etc. 
The coordinate ring of the curve ${C_0}$ determines a scalar diagonal
stabilizer algebra
\begin{equation}
\begin{aligned} 
A_0 &= \pi^*\big(H^0({C_0}\setminus\{p\},\mathcal{O}_{{C_0}})\big)\\
&\subset \pi^*\big(H^0({U_p}\setminus\{p\},\mathcal{O}_{U_p})\big)\\
&\subset H^0\big(U_o\setminus\{o\},\mathcal{O}_{U_o}\big)\\
&= L \subset gl(n,L)
\end{aligned}\label{3.3}
\end{equation}
satisfying that $A_0\cdot W\subset W$, where $L$ is identified
 with the set of scalar matrices in $gl(n, L)$.
 The rank of $W$ over $A_0$
is $r\cdot n$, which is equal to the rank of $f_*\mathcal{F}$.
Note that we have also an inclusion
$$
A_0 \isom H^0({C_0}\setminus \{p\},\mathcal{O}_{{C_0}})\subset
H^0(U_p\setminus \{p\},\mathcal{O}_{U_p}) = k((y))
$$
by the coordinate $y$. As in Section~2 and 3 of \cite{M3}, 
we can use the formal patching
$C_0 = (C_0\setminus \{p\})\cup U_p$ to compute the cohomology
group 
\begin{equation}
\begin{aligned}
H^1(C_0, \mathcal{O}_{C_0}) &\isom 
\frac{H^0(U_p\setminus \{p\}, \mathcal{O}_{U_p})}{H^0(C_0\setminus \{p\},
\mathcal{O}_{C_0}) + H^0(U_p,\mathcal{O}_{U_p})}\\
&\isom \frac{k((y))}{A_0 + k[[y]]}\;.
\end{aligned} \label{3.4}
\end{equation}
Thus the cokernel of the projection $\gam_{A_0}:A_0\longrarrow 
k((y))/k[[y]]$ has finite dimension. The function ring
$$
A_\n =H^0(C_\n\setminus \Delta, \mathcal{O}_{C_\n})
\subset \bigoplus_{j=1}^\ell  H^0(U_j\setminus
\{p_j\}, \mathcal{O}_{U_j})
$$
also acts on $V$ and satisfies that $A_\n\cdot W\subset W$,
because we have a natural injective isomorphism
\begin{equation}
\begin{aligned}
A_\n = H^0(C_\n\setminus \Delta, \mathcal{O}_{C_\n})
&\isom H^0\big(C_0\setminus \{p\}, f_*\mathcal{O}_{C_\n}\big)\\
&\subset H^0\big(U_p\setminus\{p\}, (f_*\mathcal{O}_{C_\n})_{U_p}\big)\\
&= H^0\big(U_p\setminus\{p\}, \bigoplus_{j=1}^\ell
f_{j*}\mathcal{O}_{U_j}\big)\\
&=\bigoplus_{j=1}^\ell k((y))\big[ h_{n_j}(y)\big]\\
&= H_\n(y)\subset gl\big(n,k((y))\big)\;,
\end{aligned} \label{3.5}
\end{equation}
where $h_{n_j}(y)$ is the block matrix of (\ref{2.8}) 
and $H_\n(y)$ is the maximal commutative subalgebra
of $gl\big(n,k((y))\big)$ of type $\n$.
In order to see the action of $A_\n$ on $W$ more
explicitly, we first note that the above isomorphism is given by
the identification
$y^{1/n_j} = h_{n_j}(y)$. Since
 the formal completion $\mathcal{F}_{U_j}$ of the vector bundle
$\mathcal{F}$ at the point $p_j$ is a free $\mathcal{O}_{U_j}$-module
of rank $r$, let us take a basis $\{e_1, e_2, \cdots,
e_r\}$ for the free $H^0(U_j,\mathcal{O}_{U_j})$-module  
$H^0(U_j,\mathcal{F}_{U_j})$.
The direct image sheaf $f_{j*}\mathcal{F}_{U_j}$ is a
free $\mathcal{O}_{U_p}$-module of rank $n_j\cdot r$, so we can take
a basis of sections
\begin{equation}
\left\{ y^{\alpha/n_j}\tensor e_\beta\right\}_{0\le \alpha< n_j,
1\le\beta\le r} \label{3.6}
\end{equation}
for the free $H^0(U_p,\mathcal{O}_{U_p})$-module 
$H^0(U_p,f_{j*}\mathcal{F}_{U_j})$. 
Since $H^0(U_j,\mathcal{F}_{U_j})$ $=$
$H^0(U_p,f_{j*}\mathcal{F}_{U_j})$, $H^0(U_j,\mathcal{O}_{U_j})$ $=$
$H^0(U_p,f_{j*}\mathcal{O}_{U_j})$ acts on
the basis (\ref{3.6}) by the matrix $h_{n_j}(y)\tensor I_r$,
where $I_r$ is the identity matrix acting on $\{e_1, e_2, \cdots,
e_r\}$. This can be understood by observing that the action of 
$y^{1/n}$ on the vector
$$
(c_0, c_1, \cdots, c_{n-1}) = \sum_{\alpha=0}^{n-1} 
c_\alpha y^{\alpha/n}
$$
is given by the action of the block matrix $h_n(y)$.

\medskip
\begin{rem} \label{{3.7.} Remark}
{From} the above argument, it is clear that the role
which our $\pi$ and $\phi$ play is exactly the same as that
of the \emph{parabolic structure} of \cite{Mum2}. The 
advantage of using $\pi$ and $\phi$ rather than the parabolic 
structure lies in their functoriality. Indeed, the 
parabolic structure does not transform functorially under 
morphisms of curves, while our data naturally do (see Definition~\ref{3.14. Definition}).
\end{rem}
\medskip

The algebra $H_\n(y)$ has two different presentations in
terms of geometry. We have used 
$$
H_\n(y) \isom H^0\big(U_p\setminus\{p\}, 
(f_*\mathcal{O}_{C_\n})_{U_p}\big) = \bigoplus_{j=1}^\ell 
 k((y))\big[h_{n_j}(y)\big]
\subset gl\big(n,k((y))\big)
$$
in (\ref{3.5}). In this presentation, an element of $H_\n(y)$ is an
$n\times n$ matrix acting on 
$V\isom H^0\big(U_p\setminus\{p\}, (f_*\mathcal{F})_{U_p}\big)$. The other
geometric interpretation is
$$
H_\n(y)\isom H^0\big(U_p\setminus\{p\}, \bigoplus_{j=1}^\ell
f_{j*}\mathcal{O}_{U_j}\big)\isom \bigoplus_{j=1}^\ell
H^0\big(U_j\setminus\{p_j\},\mathcal{O}_{U_j}\big)
=\bigoplus_{j=1}^\ell k((y_j))\;.
$$
In this presentation, the algebra $H_\n(y)$ acts on
$$
\begin{aligned}
V &\isom H^0\bigg(U_p\setminus\{p\},
\pi_*\bigoplus_{j=1}^\ell \psi_{j*}\big(\mathcal{O}_{U_{oj}}(-1)\big) \bigg)\\
&\isom \bigoplus_{j=1}^\ell H^0\big(U_{oj}\setminus \{o\},
\mathcal{O}_{U_{oj}}(-1)\big)\\
&= \bigoplus_{j=1}^\ell k((z_j))
\end{aligned}
$$
by the component-wise multiplication of $y_j$ to $z_j$. 
We will come back to this point in (\ref{4.4}).

The pull-back through the morphism $f$ gives an embedding
$A_0\subset A_\n$.
As an $A_0$-module, $A_\n$ is torsion-free of rank $n$, because
$C_0$ is integral and the morphism $f$ is of degree $n$. 
Using the formal patching $C_\n = (C_\n\setminus \Delta)
\cup U_1\cup\cdots\cup U_\ell$,
we can compute the cohomology
\begin{equation}
\begin{aligned}
H^1(C_\n,\mathcal{O}_{C_\n}) &\isom 
\frac{\bigoplus_{j=1}^\ell H^0(U_j\setminus \{p_j\}, 
\mathcal{O}_{U_j})}{H^0(C_\n\setminus \Delta,\mathcal{O}_{C_\n})
+ \bigoplus_{j=1}^\ell H^0(U_j,\mathcal{O}_{U_j})}\\
&\isom \frac{\bigoplus_{j=1}^\ell k((y^{1/n_j}))}{A_\n 
+ \bigoplus_{j=1}^\ell k[[y^{1/n_j}]]}\\
&\isom \frac{H_\n(y)}{A_\n + H_\n(y)\cap gl\big(n, k[[y]]\big)}\;.
\end{aligned} \label{3.8}
\end{equation}
This shows that the projection
$$
\gam_{A_\n}: A_\n\longrarrow
\frac{H_\n(y)}{H_\n(y)\cap gl\big(n, k[[y]]\big)}
$$
has a finite-dimensional cokernel. These 
discussions motivate the following definition:

\medskip
  
\begin{Def} \label{3.9. Definition}
A triple $(A_0,A_\n,W)$ is said to be a set of \emph{algebraic data}
of type $\n$, index $\mu$, and rank $r$ if the 
following conditions are satisfied:
\begin{enumerate}

\item $W$ is a point of the Grassmannian $Gr_n(\mu)$ of index $\mu$
of the vector valued functions of size $n$.

\item The type $\n$ is an integral vector $(n_1,\cdots,
n_\ell)$ consisting of positive integers  such that $n
= n_1 + \cdots + n_\ell$.

\item There is a monic element $y\in L = k((z))$ of order $-r$
such that $A_0$ is a subalgebra of $k((y))$ containing the field
$k$.

\item  The cokernel
of the projection $\gam_{A_0}:A_0\longrarrow k((y))/k[[y]]$ has 
finite dimension.

\item $A_\n$ is a subalgebra
of the maximal commutative algebra $H_\n(y) \subset
gl\big(n,k((y))\big)$ of 
type $\n$ such that the projection 
$$
\gam_{A_\n}:A_\n\longrarrow \frac{H_\n(y)}{H_\n(y)\cap
gl\big(n,k[[y]]\big)}
$$
has a finite-dimensional cokernel.

\item There is an embedding $A_0\subset A_\n$ as the scalar
diagonal matrices, and as an $A_0$-module (which is
automatically torsion-free), $A_\n$ has rank $n$ over $A_0$.

\item The algebra $A_\n\subset gl\big(n,k((y))\big)$
stabilizes $W\subset V$, i.e.\ $A_\n\cdot W\subset W$.
\end{enumerate}
\end{Def}
\medskip\noindent
 The homomorphisms $\gam_{A_0}$ and $\gam_{A_\n}$
satisfy the Fredholm condition because (7) implies that they have 
finite-dimensional kernels.
Now we can state 
\medskip
\begin{prop} \label{3.10. Proposition}
For every set of geometric data of Definition~\ref{3.1. Definition}, there is a 
unique set of algebraic data of Definition~\ref{3.9. Definition} having the same
type, index and rank.
\end{prop}
\medskip
\begin{proof}
We have already constructed the triple $(A_0,A_\n,W)$ out of
the geometric data in (\ref{3.2}), (\ref{3.3}) and (\ref{3.5}) which satisfies
all the conditions in Definition~\ref{3.9. Definition} but (1). The only remaining thing we have 
to show is that the vector subspace $W$ of (\ref{3.2}) is indeed
a point of the Grassmannian $Gr_n(\mu)$. To this end, we
need to compute the cohomology of $f_*\mathcal{F}$ 
by using the formal patching $C_0 = \Spec(A_0)\cup U_p$
(for more detail, see \cite{M3}).
Noting the identification
$$
\bigoplus_{j=1}^\ell \psi_{j*}\big(\mathcal{O}_{U_{oj}}(-1)\big) \isom 
\mathcal{O}_{U_o}(-1)^{\dsum n}
$$
as in (\ref{3.2}), we can show  that
\begin{equation}
\begin{aligned}
H^0(C_0,f_*\mathcal{F}) &= H^0(C_0\setminus \{p\},f_*\mathcal{F})
\cap H^0(U_p,f_*\mathcal{F}_{U_p})\\
&\isom W\cap H^0\big(U_p,\pi_*(\mathcal{O}_{U_o}(-1)^{\dsum n})\big)\\
&\isom W\cap H^0\big(U_o,\mathcal{O}_{U_o}(-1)^{\dsum n}\big)\\
&\isom W\cap \big(k[[z]]z\big)^{\dsum n}\\
&= \Ker(\gam_W)\;, 
\end{aligned} \label{3.11}
\end{equation}
and
\begin{equation}
\begin{aligned}
H^1(C_0,f_*\mathcal{F}) &\isom \frac{H^0(U_p\setminus \{p\}, 
f_*\mathcal{F})}{H^0(C_0\setminus \{p\},f_*\mathcal{F})
+ H^0(U_p,f_*\mathcal{F}_{U_p})}\\
&\isom \frac{H^0\big(U_p\setminus \{p\}, 
\pi_*(\mathcal{O}_{U_o}(-1)^{\dsum n})\big)}{W +
 H^0\big(U_p,\pi_*(\mathcal{O}_{U_o}(-1)^{\dsum n})\big)}\\
&\isom \frac{H^0\big(U_o\setminus \{o\}, 
\mathcal{O}_{U_o}(-1)^{\dsum n}\big)}{W +
 H^0\big(U_o,\mathcal{O}_{U_o}(-1)^{\dsum n}\big)}\\
&\isom \frac{k((z))^{\dsum n}}{W + \big(k[[z]]z\big)^{\dsum n}}\\
&=\Coker(\gam_W)\;,
\end{aligned} \label{3.12}
\end{equation}
where $\gam_W$ is the canonical projection of Definition~\ref{2.2. Definition}. 
Since $f$ is a finite morphism, we have 
$H^i(C_0,f_*\mathcal{F}) \isom H^i(C_\n, \mathcal{F})$.
Thus
\begin{equation}
\mu = \dim_k H^0(C_\n,\mathcal{F}) - \dim_k H^1(C_\n, \mathcal{F})
= \dim_k \Ker(\gam_W) - \dim_k \Coker(\gam_W)\;,\label{3.13}
\end{equation}
which shows that $W$ is indeed a point of $Gr_n(\mu)$.
This completes the proof.
\end{proof}
\medskip
\noindent
This proposition gives a generalization of the Segal-Wilson map
to the case of covering morphisms of algebraic curves. We can
make the above map further into a functor, which we shall
call the \emph{cohomology functor for covering morphisms}. The
categories we use are the following:
\medskip

\begin{Def} \label{3.14. Definition}
The category $\mathcal{C}(\n)$ of geometric data of 
a fixed type $\n$ consists of the set 
of geometric data of type $\n$ and arbitrary
 index $\mu$ and rank $r$ as its object. A morphism between two objects
$$
\left\langle f:\big( C_\n, \Delta, \Pi, \mathcal{F}, \Phi\big)
\longrarrow \big({C_0}, p, \pi, f_*\mathcal{F}, \phi\big)\right\rangle
$$
of type $\n$, index $\mu$ and rank $r$ and
$$
\left\langle {f'}:\big( C'_\n, \Delta', \Pi', \mathcal{F}', \Phi'\big)
\longrarrow \big({C'_0}, p', \pi', {f'}_*\mathcal{F}', \phi'\big)\right\rangle
$$
of the same type $\n$,  index $\mu'$ and rank $r'$
is a triple  $(\alpha,\beta, \lambda)$ of morphisms
satisfying the following conditions:
\begin{enumerate}

\item $\alpha:C'_0\longrarrow C_0$ is a locally cyclic covering
of degree $s$ of the base curves such that $\alpha^{*}(p) = 
s\cdot p'$, and  $\pi$ and $\pi'$ are related by
$\pi = \widehat\alpha\circ\pi'$
with the  morphism $\widehat\alpha$ 
of formal schemes induced by $\alpha$.

\item $\beta:C'_\n \longrarrow C_\n$ is a covering morphism of degree
$s$ such that $\Delta' = \beta^{-1}(\Delta)$,  and the 
following diagram
$$
\begin{CD}
C'_\n @>{\beta}>> C_\n\\
@V{{f'}}VV @VV{f}V\\
C'_0 @>>{\alp}> C_0
\end{CD}
$$
commutes.

\item The morphism $\widehat\beta_j:U'_j\longrarrow U_j$ of formal
schemes induced by $\beta$ at each $p'_j$ satisfies 
$\pi_j = \widehat\beta_j\circ \pi'_j$ and the 
commutativity of
$$
\begin{CD}
U_{oj} @>{\pi'_j}>> U'_j @>{\widehat\beta_j}>> U_j\\
@V{\psi_j}VV @V{{f'}_j}VV @V{f_j}VV\\
U_o @>>{\pi'}> U'_{p'} @>>{\widehat\alpha}> U_p.
\end{CD}
$$

\item $\lam: \beta_*\mathcal{F}'\longrarrow \mathcal{F}$ is an injective
$\mathcal{O}_{C_\n}$-module homomorphism such that its completion $\lam_j$
at each point $p_j$ satisfies commutativity of
$$
\begin{CD}
(\beta_*\mathcal{F}')_{U_j} @>{\lam_j}>> \mathcal{F}_{U_j}\\
@V{\widehat\beta_j(\phi'_j)}V{\wr}V @V{\wr}V{\phi_j}V\\
\widehat\beta_{j*}\pi'_{j*}\mathcal{O}_{U_{oj}}(-1) @= 
\pi_{j*}\mathcal{O}_{U_{oj}}(-1).
\end{CD}
$$
In particular, each $\lam_j$ is an isomorphism.
\end{enumerate}
\end{Def}
\medskip
\begin{rem} \label{{3.15.} Remark}
{From} (3) above, we have $r = s\cdot r'$. The condition (4) above
implies that $\mathcal{F}\big/\beta_*\mathcal{F}'$ is a torsion sheaf on
$C_\n$ whose support does not intersect with $\Delta$.
\end{rem}
\medskip\noindent
One can show by using Proposition~\ref{1.6. Proposition}
that there are many 
nontrivial morphisms among the sets of 
geometric data with different ranks.
\medskip

\begin{Def} \label{3.16. Definition}
The category $\mathcal{S}(\n)$ of algebraic data of type $\n$
has the
stabilizer triples $(A_0,A_\n,W)$ of Definition~\ref{3.9. Definition} of type $\n$
and  arbitrary  index $\mu$ and
rank $r$ as its objects. Note that for every object
$(A_0,A_\n,W)$, we have the commutative algebras
$k((y))$ and $H_\n(y)$ associated with it. A morphism between two
objects $(A_0,A_\n,W)$ and $(A'_0,A'_\n,{W'})$ is a triple
$(\iota, \epsilon, \omega)$ of injective homomorphisms
satisfying the following conditions:
\begin{enumerate}

\item $\iota:A_0\hookrarrow A'_0$ is an inclusion compatible with
the inclusion $k((y))\subset k(({y'}))$ defined by a power series
$$
y = y({y'}) = {y'}^s + a_1 {y'}^{s+1} + a_2 {y'}^{s+2} + \cdots\;.
$$

\item $\epsilon:A_\n\longrarrow A'_\n$ is an injective homomorphism
satisfying the commutativity of the diagram
$$
\begin{CD}
A_\n @>{\epsilon}>> A'_\n\\
@VVV @VVV\\
H_\n(y) @>>{\mathcal{E}}> H_\n({y'}),
\end{CD}
$$
where the vertical arrows are the inclusion maps, and 
$$
\mathcal{E}: H_\n(y) \isom \bigoplus_{j=1}^\ell k((y^{1/n_j}))
\longrarrow \bigoplus_{j=1}^\ell k(({y'}^{1/n_j})) \isom H_\n({y'})
$$
is an injective homomorphism defined by the Puiseux expansion 
$$
y^{1/n_j} = y({y'})^{1/n_j} = {y'}^{s/n_j} + b_1 {y'}^{(s+1)/n_j}
+ b_2 {y'}^{(s +2)/n_j} + \cdots
$$
of $(1)$ for every $n_j$. Note that neither 
$\epsilon$ nor $\mathcal{E}$ is an
inclusion map of subalgebras of $gl(n,L)$.

\item $\omega: {W'}\longrarrow W$ is an injective $A_\n$-module
homomorphism. We note that ${W'}$ has a natural $A_\n$-module
structure by the homomorphism $\epsilon$. As in $(2)$, $\omega$
is not an inclusion map of the vector subspaces of $V$.
\end{enumerate}
\end{Def}
\medskip
\begin{thm} \label{3.17. Theorem}
There is a fully-faithful functor 
$$
\kappa_{\n}:\mathcal{C}(\n)\overset{\sim}{\longrarrow}\mathcal{S}(\n)
$$
between the category of geometric data and the category of 
algebraic data. An object of $\mathcal{C}(\n)$ of index $\mu$ and
rank $r$ corresponds to an object of $\mathcal{S}(\n)$ of the same index
and rank.
\end{thm}
\medskip
\begin{proof}
The association of $(A_0,A_\n,W)$ to the geometric data has been done
in (\ref{3.2}), (\ref{3.3}), (\ref{3.5}) and Proposition~\ref{3.10. Proposition}. Let $(\alp,\beta,\lam)$ be a morphism
between two sets of geometric data as in Definition~\ref{3.14. Definition}. We use the
notations $U_j^* = 
U_j\setminus \{p_j\}$ and $U_p^* = U_p\setminus \{p\}$.
The homomorphism $\iota$ is defined by the commutative diagram
$$
\begin{CD}
A_0 @>{\sim}>> H^0(C_0\setminus \{p\},\mathcal{O}_{C_0}) @>>>
H^0(U_p^*,\mathcal{O}_{U_p})\\
@V{\iota}VV @V{\alp^*}VV @V{\widehat{\alp}^*}VV\\
A'_0 @>{\sim}>> H^0(C'_0\setminus \{p'\},\mathcal{O}_{C'_0}) @>>>
H^0({U'}_{p'}^*, \mathcal{O}_{{U'}_{p'}}).
\end{CD}
$$
Similarly, 
$$
\begin{CD}
A_\n @>{\sim}>> H^0(C_\n\setminus \Delta, \mathcal{O}_{C_\n}) @>>>
\bigoplus_{j=1}^\ell H^0(U_j^*,\mathcal{O}_{U_j})\\
@V{\epsilon}VV @V{\beta^*}VV @V{\oplus\widehat{\beta}_j^*}VV\\
A'_\n @>{\sim}>> H^0(C'_\n\setminus \Delta',\mathcal{O}_{C'_\n}) @>>>
\bigoplus_{j=1}^\ell H^0({U'}_{j}^*, \mathcal{O}_{{U'}_{j}})
\end{CD}
$$
defines the homomorphism $\epsilon$. Finally, 
$$
\begin{CD}
{W'} @= {W'} @>{\omega}>> W\\
@VV{\wr}V @VV{\wr}V @VV{\wr}V\\
H^0(C'_0\setminus \{p'\}, {{f'}}_*\mathcal{F}') @>{\alp_*}>{\sim}>
H^0(C_0\setminus \{p\}, f_*\beta_*\mathcal{F}') @>{f_*(\lam)}>>
H^0(C_0\setminus \{p\}, f_*\mathcal{F})\\
@V{{f'}^*}V{\wr}V @V{f^*}V{\wr}V @V{f^*}V{\wr}V\\
H^0(C'_\n\setminus \Delta', \mathcal{F}') @>{\beta_*}>{\sim}>
H^0(C_\n\setminus \Delta, \beta_*\mathcal{F}') @>{\lam}>>
H^0(C_\n\setminus \Delta, \mathcal{F})\\
@VVV @VVV @VVV\\
\bigoplus_j H^0({U'}_j^*, \mathcal{F}'_{U'_j}) @>{\dsum
\widehat{\beta}_{j*}}>{\sim}>
\bigoplus_j H^0(U_j^*, \widehat{\beta}_{j*}\mathcal{F}'_{U'_j}) 
@>{\dsum\lam_j}>>
\bigoplus_j H^0(U_j^*, \mathcal{F}_{U_j})
\end{CD}
$$
determines the homomorphism $\omega$. 

In order to establish that the two categories are equivalent,
we need the inverse construction. The next section is entirely devoted
to the proof of this claim.
\end{proof}
\medskip

 The following proposition and its corollary about the geometric
 data of rank one are  crucial when we study geometry of orbits of the
 Heisenberg flows in Section~\ref{sec: Prym}.
\medskip
\begin{prop} \label{3.18. Proposition}
Suppose we have two sets of geometric data of rank one having exactly
the same constituents except for the sheaf isomorphisms 
$(\Phi ,\phi)$ for one and $(\Phi' ,\phi')$ for the other. Let
$(A_0,A_\n,W)$ and $(A_0,A_\n,{W'})$ be the corresponding
algebraic data, where $A_0$ and $A_\n$ are common in both of
the triples because of the assumption. Then there is an
element $g\in \Gam_\n(y)$ of $(\ref{2.15})$ such that ${W'} = g\cdot W$.
\end{prop}
\medskip
\begin{proof}
Recall that 
$$
\phi:(f_*\mathcal{F})_{U_p}\overset{\sim}{\longrarrow}
\pi_*\bigg(\bigoplus_{j=1}^\ell \psi_{j*}\big(\mathcal{O}_{U_{oj}}(-1)
\big)\bigg)
$$
is an $\big(f_*\mathcal{O}_{C_\n}\big)_{U_p}$-module isomorphism.
  Thus,
$$
g = \phi'\circ \phi^{-1}: \pi_*\bigg(\bigoplus_{j=1}^\ell 
\psi_{j*}\big(\mathcal{O}_{U_{oj}}(-1)\big)\bigg)
\overset{\sim}{\longrarrow}
\pi_*\bigg(\bigoplus_{j=1}^\ell \psi_{j*}
\big(\mathcal{O}_{U_{oj}}(-1)\big)\bigg)
$$
is also an
 $\big(f_*\mathcal{O}_{C_\n}\big)_{U_p}$-module isomorphism. 
Note that 
we have identified  $\big(f_*\mathcal{O}_{C_\n}\big)_{U_p}$ 
as a subalgebra of $H_\n(y)$ in (\ref{3.5}). 
Indeed, this subalgebra is $H_\n(y)\cap
gl\big(n,k[[y]]\big)$. Therefore, the invertible $n\times n$ matrix
$$
g \in k^{*\dsum n} + gl\big(n,k[[y]]y\big)
= k^{*\dsum n}  + gl\big(n,k[[z]]z\big)
$$
 commutes with $H_\n(y)\cap
gl\big(n,k[[y]]\big)$, where $k^*$ denotes the set of nonzero
constants and $k^{*\dsum n}$ the set of invertible constant
diagonal matrices. We recall that $k[[z]] = k[[y]]$, because $y$ has
order $-1$. The commutativity of $g$ and $H_\n(y)\cap
gl\big(n,k[[y]]\big)$ immediately implies that $g$ commutes
with all of $H_\n(y)$. But since $H_\n(y)$ is a \emph{maximal}
commutative subalgebra of $gl\big(n,k((y))\big)$, it implies
that $g\in \Gam_\n(y)$. Here we note that $\phi'_j\circ {\phi_j}^{-1}$
is exactly the $j$-th block of size
$n_j\times n_j$ of the $n\times n$ matrix $g$, and that 
we can normalize the leading term of  $\phi'_j\circ {\phi_j}^{-1}$
to be equal to $I_{n_j}$ by the definition
(5) of Definition~\ref{3.1. Definition}. Thus the leading term of $g$ can be normalized
to $I_n$.  {From} the construction of (\ref{3.2}), we have
${W'} = g\cdot W$. This completes the proof.
\end{proof}
\medskip
\begin{cor} \label{3.19. Corollary}
The cohomology functor induces a bijective correspondence
 between the collection
of geometric data
$$
\left\langle f:\big( C_\n, \Delta, \Pi, \mathcal{F}\big)
\longrarrow \big({C_0}, p, \pi, f_*\mathcal{F}\big)\right\rangle
$$
of type $\n$, index $\mu$, and rank one, and the triple
of algebraic data $(A_0,A_\n,\overline{W})$ of type $\n$,
index $\mu$, and rank one 
satisfying the same conditions of Definition~\ref{3.9. Definition} except that
 $\overline{W}$ is a point of the Grassmannian quotient
$Z_\n(\mu,y)$. 
\end{cor}
\medskip
\begin{proof}
Note that the datum $\Phi$ is indeed the block decomposition
of the datum of $\phi$. Thus taking the quotient space of the
Grassmannian by the group action of $\Gam_n(y)$
exactly corresponds to eliminating the data
$\Phi$ and $\phi$ from the set of geometric data of Definition~\ref{3.1. Definition}.
\end{proof}

\bigskip

\section{The inverse construction}
\label{sec: inverse}
\medskip

Let $W\in Gr_n(\mu)$ be a point of the Grassmannian and
consider a commutative subalgebra $A$ of $gl(n,L)$
such that $A\cdot W\subset W$. Since the set of vector fields
${\Psi}(A)$ has $W$ as a fixed point, we call such an algebra a
commutative \emph{stabilizer} algebra of $W$. In the previous
work \cite{M3}, the algebro-geometric structures
of arbitrary commutative stabilizers
were determined for the case of
the Grassmannian $Gr_1(\mu)$ of scalar valued
functions. In the context of the current paper, the Grassmannian is
enlarged, and consequently there are far larger varieties
of commutative stabilizers. However, it is not the  purpose 
of this paper to give the
complete geometric classification
of arbitrary stabilizers. We restrict ourselves to
studying \emph{large} stabilizers in connection with Prym 
varieties, which will be the central theme of the next section. A
stabilizer is said to be large if it corresponds to a finite-dimensional
orbit of the Heisenberg flows on the Grassmannian quotient.
The goal of this section is to recover the geometric data
out of a point of the Grassmannian together with a large stabilizer.

 Choose an integral vector
$\n = (n_1, n_2, \cdots, n_\ell)$ with  $n = n_1 + \cdots + n_\ell$
and a monic element $y$ of order $-r$ as in (\ref{2.7}), and
consider the formal loop algebra $gl\big(n,k((y))\big)$ acting
on the vector space $V = L^{\dsum n}$. 
Let us denote $y_j = h_{n_j}(y) = y^{1/n_j}$. We introduce a 
filtration 
$$
\cdots \subset H_\n(y)^{(rm-r)}\subset H_\n(y)^{(rm)}\subset
H_\n(y)^{(rm+r)}\subset \cdots
$$
in the maximal commutative algebra
\begin{equation}
H_\n(y) \isom \bigoplus_{j=1}^\ell k((y))\big[y^{1/n_j}\big]
\isom \bigoplus_{j=1}^\ell k((y^{1/n_j}))
= \bigoplus_{j=1}^\ell k((y_j))\label{4.1}
\end{equation}
by defining 
\begin{equation}
H_\n(y)^{(rm)} = \bigg\{\big(a_1(y_1), \cdots, a_{\ell}(y_{\ell})\big)\;
\bigg| \;\max\big[\ord_{y_1}(a_1),\cdots ,\ord_{y_\ell}(a_\ell)\big]
\le m\bigg\}\;,\label{4.2}
\end{equation}
where $\ord_{y_j}(a_j)$ is the order of $a_j(y_j)\in k((y_j))$ with
respect to the variable $y_j$. Accordingly, we can introduce
a filtration in $V$ which is compatible with the
action of $H_\n(y)$ on $V$.  In order to define the new filtration in 
$V$ geometrically, let us start with $U_o = \Spec\big(k[[z]]\big)$
and $U_p = \Spec\big(k[[y]]\big)$. The inclusion $k[[y]]\subset
k[[z]]$ given by $y = y(z) = z^r$ $+$ $c_1z^{r+1}$ $+$ $c_2z^{r+2}$
$+ \cdots$ defines a morphism $\pi:U_o\longrarrow U_p$.
Let $U_j = \Spec\big(k[[y_j]]\big)$. The identification $y_j =
y^{1/n_j}$ gives a cyclic covering $f_j:U_j\longrarrow U_p$
of degree $n_j$. Correspondingly,
the covering $\psi_j:U_{oj}\longrarrow U_o$ of
degree $n_j$ of (9) of Definition~\ref{3.1. Definition} is given by $k[[z]]\subset k[[z^{1/n_j}]]$.
Thus we have a commutative diagram
$$
\begin{CD}
k[[z^{1/n_j}]] @<{\pi^*_j}<< k[[y^{1/n_j}]]\\
@A{\psi^*_j}AA @AA{f^*_j}A\\
k[[z]] @<<{\pi^*}< k[[y]]
\end{CD}
$$
of inclusions, where $\pi^*_j$ is defined by the Puiseux expansion
\begin{equation}
y_j = y^{1/n_j} = y(z)^{1/n_j} = z^{r/n_j} + a_1z^{(r+1)/n_j}
+ a_2z^{(r+2)/n_j} + \cdots\label{4.3}
\end{equation}
of $y(z)$. Recall that 
in order to distinguish from  $U_o = \Spec\big(k[[z]]\big)$,
we have introduced the notation
$U_{oj}$ $=$ $\Spec\big(k[[z^{1/n_j}]]\big)$ for the 
cyclic covering of $U_o$.
The above diagram corresponds to the geometric diagram of covering
morphisms
$$
\begin{CD}
U_{oj} @>{\pi_j}>> U_j\\
@V{\psi_j}VV @VV{f_j}V\\
U_o @>>{\pi}> U_p.
\end{CD}
$$
We denote $U^*_o = U_o\setminus \{o\}$, $U_{oj}^*
= U_{oj}\setminus \{o\}$, $U_p^*=
U_p\setminus \{p\}$, and $U_j^* = U_j\setminus \{p_j\}$ as before.
The $k((y))$-algebra $H_\n(y)$ is identified with the 
$H^0(U_p^*,\mathcal{O}_{U_p})$-algebra 
$$
H_\n(y) = H^0\big(U_p^*, \bigoplus_{j=1}^\ell f_{j*}\mathcal{O}_{U_j}\big)
\isom \bigoplus_{j=1}^\ell H^0(U_j^*,\mathcal{O}_{U_j})\;.
$$
Corresponding to this identification, the vector space 
$V = L^{\dsum n}$ as a module over $L = H^0(U_o^*,\mathcal{O}_{U_o})$
is identified with
\begin{equation}
V = H^0\bigg(U_o^*,\bigoplus_{j=1}^\ell 
\psi_{j*}\big(\mathcal{O}_{U_{oj}}(-1)\big)\bigg) \isom \bigoplus_{j=1}^\ell
H^0\big(U_{oj}^*,\mathcal{O}_{U_{oj}}(-1)\big) 
\isom \bigoplus_{j=1}^\ell k((z^{1/n_j}))\;.
\label{4.4}
\end{equation}
The $H_\n(y)$-module structure of $V$ is given by the pull-back
$\bigoplus_{j=1}^\ell \pi_j^*$, which is nothing but the 
component-wise multiplication of $k((y^{1/n_j}))$ to 
$k((z^{1/n_j}))$ through (\ref{4.3}) for each $j$. 
Define a new variable by $z_j = z^{1/n_j}$.
We note from (\ref{4.3}) that $y_j= y_j(z^{1/n_j}) = y_j(z_j)$
is of order $-r$ with respect to $z_j$. Now we can introduce 
a new filtration 
$$
\cdots \subset V^{(m-1)}\subset V^{(m)}\subset V^{(m+1)}\subset \cdots
$$ 
in $V$ by defining
\begin{equation}
V^{(m)} = \bigg\{\big(v_1(z_1), \cdots, v_{\ell}(z_{\ell})\big)
\in \bigoplus_{j=1}^\ell k((z_j))\;
\bigg| \;
\max\big[\ord_{z_1}(v_1),\cdots ,\ord_{z_\ell}(v_\ell)\big]
\le m\bigg\}\;,\label{4.5}
\end{equation}
where $\ord_{z_j}(v_j)$ denotes the order of $v_j = v_j(z_j)$
with respect to $z_j$.
\medskip
\begin{rem} \label{{4.6.} Remark}
The filtration (\ref{4.5}) is different from (\ref{2.1}) in general. However, 
we always have  $V^{(0)} = F^{(0)}(V)$ and $V^{(-1)} = 
F^{(-1)}(V)$. This is one of the reasons why we have chosen 
$F^{(-1)}(V)$ instead of an arbitrary $F^{(\nu)}(V)$
in the definition of the Grassmannian in Definition~\ref{2.2. Definition}.
\end{rem}
\medskip
 It is clear from (\ref{4.2}) and (\ref{4.5}) that  
$H_\n(y)^{(rm_1)}\cdot V^{(m_2)}
\subset V^{(rm_1+m_2)}$,
and hence $V$ is a filtered $H_\n(y)$-module. With these preparation, 
we can state the inverse construction theorem.
\medskip
\begin{thm} \label{4.7. Theorem}
A triple $(A_0,A_\n,W)$ of algebraic data of Definition~\ref{3.9. Definition} determines
a unique set of geometric data
$$
\left\langle f:\big( C_\n, \Delta, \Pi, \mathcal{F}, \Phi\big)
\longrarrow \big({C_0}, p, \pi, f_*\mathcal{F}, \phi\big)\right\rangle\;.
$$
\end{thm}
\medskip
\begin{proof} The proof is divided into four parts.
\smallskip
\noindent
(I) {\sl Construction of the curve $C_0$ and the point $p$:}
  Let us define
$A_0^{(rm)} = A_0\cap k[[y]]y^{-m}$, 
which   consists of elements of $A_0$ of order at most $m$ with 
respect to the variable $y$.
 This gives a filtration of $A_{0}$:
$$
\cdots \subset A_{0}^{(rm-r)}\subset A_{0}^{(rm)}\subset 
A_{0}^{(rm+r)}
\subset \cdots\;.
$$
 Using the  finite-dimensionality of the cokernel (4) of Definition~\ref{3.9. Definition},
 we can show  that $A_0$ has an element of order $m$ (with
respect to $y$) for
every large integer $m\in \mathbb{N}$, i.e.\
\begin{equation}
\dim_k A_0^{(rm)}\big/A_0^{(rm-r)} = 1 \quad{\text{ for all }}
\quad m>>0\;.
\label{4.8}
\end{equation}
 Since $A_0\cdot W\subset W$, the Fredholm condition of $W$ 
implies  that
$A_0^{(rm)} = 0$ for all $m<0$.
Note that $A_0$ is a subalgebra of a field, 
and hence it is an integral domain. Therefore, 
the complete algebraic curve
${C_0} = \Proj(gr A_0)$ defined by the graded algebra
$$
gr A_0=\bigoplus_{m=0}^\infty A_0^{(rm)}
$$
is integral. 
We claim that $C_0$ is a one-point completion of the affine 
curve $\Spec(A_0)$. In order to prove the claim, 
let $w$ denote the homogeneous
 element of degree one given by the image of the element 
 $1\in A_0^{(0)}$ under the inclusion $A_0^{(0)}\subset A_0^{(r)}$.
 Then the homogeneous localization $(grA_0)_{((w))}$ is 
isomorphic to 
$A_{0}$. Thus the principal open subset $D^{+}(w)$
 defined by the element $w$ is isomorphic to the affine
 curve $\Spec(A_0)$. The complement of $\Spec(A_0)$ in $C_0$ is
 the closed subset defined by $(w)$, which is nothing but
the projective scheme  
$$
\Proj\left(\bigoplus_{m=0}^\infty A_0^{(rm)}\big/A_0^{(rm-r)}\right)
$$
 given  by the associated graded algebra of $grA_0$.
Take a monic element 
 $a_m\in A_0^{(rm)}\setminus A_0^{(rm-r)}$ for every $m >>0$,
whose existence is assured by (\ref{4.8}).
Since $a_i\cdot a_j \equiv a_{i+j}$ mod $A_0^{(ri+rj-r)}$, the map
$$
 \zeta: \bigoplus _{m=0}^\infty A_0^{(rm)}/A_0^{(rm-r)}
\longrightarrow k[x]\;,
$$
which assigns $x^m$ to each $a_m$ for $m>>0$ and 0 otherwise,
 is a well-defined homomorphism of graded
 rings, where $x$ is an indeterminate.
  In fact, $\zeta$ is an isomorphism in large degrees, and hence 
we have
$$
\Proj\left(\bigoplus_{m=0}^\infty A_0^{(rm)}\big/A_0^{(rm-r)}\right)
    \isom \Proj \big(k[x]\big) = p\;.
$$
This proves the claim.

 Next we want to show that the added
 point $p$ is a smooth rational point of $C_0$. To this end, 
it is sufficient to show that
 the formal completion of the structure sheaf 
 of $C_0$ along $p$ is isomorphic to a formal
 power series ring.  Let us
 consider $(grA_0)/({w}^n)$.  The degree $m$ homogeneous piece 
 of this ring is given by 
 $A_0^{(rm)}/(w^nA_0^{(rm-rn)})$, which is
 isomorphic to $ k\cdot a_m\oplus k\cdot a_{m-1}w
 \oplus \cdots \oplus k\cdot a_{m-n+1}w^{n-1}$ 
for all $m>n>>0$. {From} this we conclude that 
$$
gr(A_0)/(w^n)\cong k[x, w]/(w^n)
$$
 in large degrees for $n>>0$. 
Therefore, taking the homogeneous localization 
 at the ideal $(w)$, we have
$$
\big(gr(A_0)/(w^n)\big)_{((w))}\cong 
k[w/x]\big/\big((w/x)^n\big)
$$
for $n>>0$. Letting $n\rightarrow \infty$ and taking the inverse
 limit of this inverse 
 system, we  see that the formal completion of
  the structure sheaf of $C_0$ along  $p$ is indeed isomorphic to
 the formal power series ring $k[[w/x]]$. 
  We can also present an affine local neighborhood 
 of the point $p$. Let $a=a(y)\in A_0$ be a monic, nonconstant
 element with the lowest order. It is unique up to the addition
 of a constant: $a(y)\mapsto a(y)+c$.  This element defines a 
 principal open subset $D^{+}(a)$ corresponding to the ring
\begin{equation}
\begin{aligned}
(grA_0)_{(a)} &= gr A_0 \big[a^{-1}\big]_0\\
&=\{a^{-i}b\;|\;b\in A_0,\;
i\ge 0,\; \ord_y(b) - i\cdot\ord_y(a) \le 0\}\\
&\subset k[[y]]\;.
\end{aligned}\label{4.9}
\end{equation}
Since the formal completion of $C_0$ along $p$ coincides with
that of $D^+(a)$ at $p$, and since the structure sheaf of the
latter is $k[[y]]$ by (\ref{4.9}), we have obtained that
$k[[w /x]]= k[[y]]$. Thus $y$ is indeed
 a formal  parameter of the curve $C_0$ at $p$. 

\smallskip
\noindent
(II) {\sl Construction of $C_\n$ and $\Delta$:}
Since $A_\n\subset H_\n(y)$, it has a filtration
$A_\n^{(rm)} = A_\n\cap H_\n(y)^{(rm)}$ induced by (\ref{4.2}). 
The Fredholm condition of $W$ again implies that
$A_\n^{(rm)} = 0$ for all $m<0$.
So let us define $C_\n = \Proj(gr A_\n)$, where
$$
gr A_\n =  \bigoplus_{m=0}^\infty A_\n^{(rm)}\;.
$$
This is a complete algebraic curve and has an  affine part  
$\Spec(A_\n)$. The complement
$C_\n\setminus \Spec(A_\n)$ is given by the projective scheme
$$
\Proj\left(\bigoplus_{m=0}^\infty 
A_\n^{(rm)}\big/A_\n^{(rm-r)}\right)\;.
$$
The finite-dimensionality (5) of Definition~\ref{3.9. Definition}
implies that for every $\ell$-tuple $(\nu_1,\cdots,\nu_\ell)$
of positive integers satisfying that $\nu_j>>0$,
the stabilizer algebra
$A_\n$ has an element of the form
$$
\big(a_1(y_1),\;\cdots,\;a_\ell(y_\ell)\big)
 \in A_\n \subset \bigoplus_{j=1}^\ell k((y_j))
$$
such that the order of $a_j(y_j)$ with respect to $y_j$ is
equal to $\nu_j$ for all $j = 1, \cdots,\ell$.
Thus for all
 sufficiently large integer $m\in\mathbb{N}$, we have an 
isomorphism
$$
A_\n^{(rm)}\big/A_\n^{(rm-r)} \isom k^{\dsum\ell}\;.
$$
Actually, by choosing a basis of $A_\n^{(rm)}\big/A_\n^{(rm-r)}$
 for each $m>>0$, 
we can  prove in the similar way as in the scalar case
 that the associated graded algebra 
$\bigoplus _{m=0}^{\infty}A_\n^{(rm)}/A_\n^{(rm-r)}$
 is isomorphic to the graded algebra 
$\bigoplus _{j=0}^{\ell} k[x_j]$ in sufficiently large degrees, 
where $x_j$'s are independent
 variables. The projective scheme of the latter graded algebra is
 an $\ell$-point scheme. 
Therefore, the curve $C_\n$ is an $\ell$-point completion of the
affine curve $\Spec(A_\n)$. Let
$$
\Delta = \{p_1,p_2,\cdots,p_\ell\} = 
\Proj\left(\bigoplus_{m=0}^\infty 
A_\n^{(rm)}\big/A_\n^{(rm-r)}\right)\;.
$$
We have to show that these points are smooth and rational.
To this end,  we investigate the completion of 
$C_\n$ along the subscheme
 $\{ p_1, p_2, \cdots , p_{\ell}\}$. Let $u$ be 
the homogeneous element of degree one
 in $A_\n^{(r)}$ given by the image of $1\in A_\n^{(0)}$ under
 the inclusion map $A_\n^{(0)}\subset A_\n^{(r)}$. Then the 
 closed subscheme (the added points) is exactly the one defined 
 by the principal homogeneous ideal $(u)$. We can prove,
 in a similar way as in (I), that 
$$
gr(A_\n)/(u^n)\cong \big(\bigoplus_{j=1}^{\ell}k[x_j]\big)[u]\big/(u^n)
 \cong \bigoplus_{j=1}^{\ell} \big(k[x_j, u_j]/(u_j^n)\big)
$$
 in large degrees for $n>>0$, where
  $x_j$'s and  $u_j$'s are independent  variables.  Letting
 $n\rightarrow \infty$ and taking the inverse limit, we conclude that
 the formal completion of the structure sheaf of
  $C_\n$ along the subscheme $\{ p_1, p_2,
 \cdots, p_{\ell}\}$ 
 is isomorphic to the direct sum $\bigoplus_{j=1}^{\ell}k[[u_j/x_j]]$.
 Thus all of these $\ell$ points are smooth and rational.
By considering the adic-completion of the ring
$$
(A_\n)_p =\big\{a^{-i} h\;|\;h\in {A_\n}^{(rm)},\;i\ge 0,\; 
m-i\cdot \ord_y(a) \le 0\big\}\;,
$$
where $a$ is as in (\ref{4.9}), we can   show that $k[[u_j/x_j]] 
= k[[y_j]]$.   So $y_j$ can be viewed as a  
formal parameter of $C_\n$ around
  the point $p_j$.

\smallskip
\noindent
(III) {\sl Construction of the morphism $f$:}
The inclusion map $A_0\hookrarrow A_\n$ gives rise to an inclusion
\begin{equation}
\bigoplus_{q=0}^\infty A_0^{(rq)}
\subset \bigoplus_{m=0}^\infty A_\n^{(rm)}\;,\label{4.10}
\end{equation}
because we have $A_0^{(rq)} \subset A_\n^{(rm)}$
for all $m\ge q\cdot \max[n_1,\cdots,n_\ell]$.
It defines a finite surjective morphism $f:C_\n\longrarrow {C_0}$.
  Using the formal parameter $y_j$, we know that the morphism
$f_j:U_j\longrarrow U_p$ of the formal completion $U_j$ of 
$C_\n$ along $p_j$ induced by $f:C_\n\longrarrow C_o$ is indeed
the cyclic covering morphism defined by $y = y_j^{n_j}$.
Since $H_\n(y)$ is a free $k((y))$-module of dimension $n$
and since the algebras $A_0$ and $A_\n$ satisfy the Fredholm 
condition described in (4), (5) and (7) of Definition~\ref{3.9. Definition},
 $A_\n$ is a torsion-free
module of rank $n$ over $A_0$. Thus the morphism $f$ has degree $n$.

\smallskip
\noindent
(IV) {\sl Construction of the sheaf $\mathcal{F}$:}
We introduce a filtration in $W\subset V$ induced by (\ref{4.5}).
The $A_\n$-module structure of 
$W$ is compatible with the 
$H_\n(y) = \bigoplus_{j=1}^\ell k((y_j))$-action on 
$V 
= \bigoplus_{j=1}^\ell k((z_j))$.
Note that we have $A_\n^{(rm_1)}\cdot W^{(m_2)}\subset W^{(rm_1+m_2)}$,
and hence $\bigoplus_{m=-\infty}^\infty W^{(m)}$ is a 
graded module over $gr A_\n$.
Let $\mathcal{F}$ be the sheaf corresponding to the shifted graded module
 $\big(\bigoplus_{m=-\infty}^{\infty}W^{(m)}\big)(-1)$,
where this shifting by $-1$ comes from our convention of Definition~\ref{2.2. Definition}.
This sheaf is an extension of the sheaf $W^\sim$ defined on the
affine curve $\Spec(A_\n)$. The graded module 
$\big(\bigoplus_{m=-\infty}^{\infty}W^{(m)}\big)(-1)$
is also a graded module over 
$gr A_0$ by
(\ref{4.10}). It gives rise to a torsion-free sheaf on ${C_0}$, which
is nothing but $f_*\mathcal{F}$. 
Let us define
$$
W_p =\big\{a^{-i}w\;|\;w\in W^{(m)},\;i\ge 0,\; 
m-i\cdot r\cdot \ord_y(a) \le -1\big\}\;,
$$
where $a$ is as in (\ref{4.9}).
Then $W_p$ is an $(A_0)_p$-module of rank $r\cdot n = r\sum n_j$.
The formal completion $(f_*\mathcal{F})_{U_p}$ of $f_*\mathcal{F}$ at the point $p$
is given by the $k[[y]]$-module $W_p\tensor_{(A_0)_p} k[[y]]$,
and the isomorphism
\begin{equation}
W_p\tensor_{(A_0)_p} k[[y]]\isom
\bigoplus_{j=1}^\ell k[[z_j]]z_j \label{4.11}
\end{equation}
gives rise to the sheaf isomorphism
$$
\phi:(f_*\mathcal{F})_{U_p}\overset{\sim}{\longrarrow}
\pi_*\bigoplus_{j=1}^\ell \psi_{j*}\big(\mathcal{O}_{U_{oj}}(-1)\big)
$$
and its diagonal blocks $\Phi = (\phi_1,\cdots,\phi_\ell)$.
Since $f_*\mathcal{F}$ has rank $r\cdot n$ over $\mathcal{O}_{C_0}$ from (\ref{4.11})
and $A_\n$ has rank $n$
over $A_0$, the sheaf $\mathcal{F}$ on $C_\n$ must have rank $r$.
The cohomology calculation of (\ref{3.11}), (\ref{3.12}) and (\ref{3.13}) shows that the
Euler characteristic of $\mathcal{F}$ is equal to $\mu$.
Thus we have constructed all of the ingredients
of the geometric data of type $\n$,
index $\mu$, and rank $r$. This completes the proof of Theorem~\ref{4.7. Theorem}.
\end{proof}
\medskip

In order to complete the proof of the categorical equivalence
of Theorem~\ref{3.17. Theorem}, we have to construct a triple $(\alpha,\beta,\lam)$
out of the homomorphisms $\iota:A_0\hookrarrow A'_0$,
$\epsilon:A_\n \longrarrow A'_\n$, and $\omega:{W'}\longrarrow W$. 
Let $s$ be the rank of $A'_0$ as an $A_0$-module. The injection
$\iota$ is associated with the inclusion $k((y))\subset k(({y'}))$, 
and the coordinate $y$ has order $-s$ with respect to ${y'}$.
Therefore, we have $r = s\cdot r'$.
Recall that the filtration we have introduced in $A_0$ is
defined by the order with respect to $y$. The homomorphism $\iota$
induces an injective homomorphism
$$
gr A_0 = \bigoplus_{m=0}^\infty A_0^{(rm)}\longrarrow
\bigoplus_{m=0}^\infty {A'}_0^{(s\cdot r'm)}\subset
\bigoplus_{m=0}^\infty {A'}_0^{(r'm)} = gr {A'}_0\;,
$$
which then defines a morphism $\alp:{C'}_0\longrarrow C_0$.

Note that the homomorphism $\epsilon$ comes from the
inclusion $k((y_j))\subset k(({y'}_j))$ for every $j$. By the 
Puiseux expansion, we see that every $y_j = y^{1/n_j}$ has order $-s$
as an element of $k(({y'}_j)) = k(({y'}^{1/n_j}))$. Thus we have
$$
gr A_\n = \bigoplus_{m=0}^\infty A_\n^{(rm)}\longrarrow
\bigoplus_{m=0}^\infty {A'}_\n^{(s\cdot r'm)}\subset
\bigoplus_{m=0}^\infty {A'}_\n^{(r'm)} = gr {A'}_\n\;,
$$
and this homomorphism defines $\beta:{C'}_\n\longrarrow C_\n$.

Finally, the homomorphism $\lam$ can be constructed as follows.
Note that $\omega$ gives an inclusion
${W'}^{(m)}\subset W^{(m)}$ as subspaces of 
$\bigoplus_{j=1}^\ell k((z_j))$ for every 
$m\in \mathbb{Z}$. Thus we have an inclusion map
$$
\bigoplus_{m=-\infty}^\infty {W'}^{(m)} \subset
\bigoplus_{m=-\infty}^\infty W^{(m)} \;,
$$
which is clearly a $gr A_\n$-module homomorphism. Thus it
induces an injective homomorphism $\lam:\beta_*\mathcal{F}'
\longrarrow \mathcal{F}$. 

One can check that the construction we have given in Section~\ref{sec: inverse}
is indeed the inverse of the map we defined in Section~\ref{sec: coh}.
Thus we have completed the entire proof of the categorical
equivalence Theorem~\ref{3.17. Theorem}.

\bigskip

\section{A characterization of arbitrary Prym varieties}
\label{sec: Prym}
\medskip

In this section, we study the geometry of finite type orbits
of the Heisenberg flows, and establish a simple characterization
theorem of arbitrary Prym varieties.
Consider the Heisenberg flows associated with $H_\n(y)$
on the Grassmannian quotient $Z_\n(\mu,y)$ and assume that 
the flows produce a finite-dimensional
orbit at a point $\overline{W}\in Z_\n(\mu,y)$. 
Then this situation corresponds to the
geometric data of Definition~\ref{3.1. Definition}:
\medskip
\begin{prop} \label{5.1. Proposition}
Let $W\in Gr_n(\mu)$ be a point of the Grassmannian at which
the Heisenberg flows of type $\n$ and rank $r$
associated with $H_\n(y)$ 
generate an orbit of finite type. Then $W$ gives rise to
a set of geometric data
$$
\left\langle f:\big( C_\n, \Delta, \Pi, \mathcal{F}, \Phi\big)
\longrarrow \big({C_0}, p, \pi, f_*\mathcal{F}, \phi\big)\right\rangle
$$
of type $\n$, index $\mu$, and rank $r$. 
\end{prop}
\medskip
\begin{proof}
Let $X_\n$ be the orbit
of the Heisenberg flows starting at $W$, and consider
the $r$-reduced KP flows associated with $k((y))$.
The finite-dimensionality of $\overline{X}_\n
= Q_{\n,y}(X_\n)$ implies that the $r$-reduced KP flows also
produce a finite type orbit $X_0$ at $W$. 
Let $A_0 = \{a\in k((y))\;|\; a\cdot W\subset W\}$ and 
$A_\n = \{ h\in H_\n(y)\;|\; h\cdot W\subset W\}$
be the stabilizer subalgebras, which satisfy
$A_0\subset A_\n$. {From} the definition 
of the vector fields Definition~\ref{2.5. Definition}, an element of $k((y))$ gives 
the zero tangent vector at $W$ if and only if it is in
$A_0$. Similarly, for an element $b\in H_\n(y)$,
${\Psi}_W(b) = 0$ if and only if $b\in A_\n$. Thus 
the tangent spaces of these orbits are given by
$$
T_{W}X_0\isom k((y))\big/A_0
\quad{\text{  and   }}\quad
T_{W}X_\n \isom H_\n(y)\big/A_\n\;.
$$
Therefore, going down to the Grassmannian quotient, 
the tangent spaces of $\overline{X}_\n$ and $\overline{X}_0
= Q_{\n,y}(X_0)$ are now given by
$$
T_{\overline{W}}\overline{X}_0 \isom 
\frac{k((y))}{A_0 + k((y))\cap gl\big(n, k[[y]]y\big)}
= \frac{k((y))}{A_0\dsum k[[y]]y}
$$
and 
$$
T_{\overline{W}}\overline{X}_\n \isom 
\frac{H_\n(y)}{A_\n + H_\n(y)\cap gl\big(n, k[[y]]y\big)} = 
\frac{H_\n(y)}{A_\n \dsum H_\n(y)^-}\;,
$$
where $\overline{W} = Q_{\n,y}(W)$, and $H_\n(y)^-$ is defined
in (\ref{2.14}).
Since both of the above sets are finite-dimensional, 
the triple $(A_0,A_\n,W)$ satisfies the cokernel conditions
(4) and (5) of Definition~\ref{3.9. Definition}. The rank condition (6) of Definition~\ref{3.9. Definition} is a 
consequence of the fact that $H_\n(y)$ has
 dimension $n$ over $k((y))$.
Therefore, applying the inverse construction of 
the cohomology functor to the triple, 
we obtain a set of geometric data. This completes the proof.
\end{proof}
\medskip
\noindent
Since $k\subset A_0\subset A_\n$, from (\ref{3.4}) and (\ref{3.8}) we obtain
\begin{equation}
 T_{\overline{W}}\overline{X}_0 \isom\frac{k((y))}{A_0\dsum k[[y]]y}
 \isom H^1({C_0},\mathcal{O}_{{C_0}})
\label{5.2}
\end{equation}
and 
\begin{equation}
T_{\overline{W}}\overline{X}_\n \isom 
\frac{H_\n(y)}{A_\n \dsum H_\n(y)^-}
\isom H^1(C_\n,\mathcal{O}_{C_\n})\;.\label{5.3}
\end{equation}
Thus we know that the genera of ${C_0}$ and $C_\n$ are equal to 
the dimension of the orbits $\overline{X}_0$ and
 $\overline{X}_\n$ on the Grassmannian quotient, respectively.
However, we cannot conclude that these orbits are actually Jacobian
varieties. The difference of the orbits and the Jacobians lies in
the deformation of the data $(\Phi,\phi)$. In order to give a 
surjective map 
from the Jacobians to these orbits, we have to eliminate these
unwanted information by using Corollary~\ref{3.19. Corollary}. Therefore, in the rest
of this section, we have to assume that the point $W\in Gr_n(\mu)$ 
gives rise to
a rank one triple $(A_0,A_\n,W)$ of algebraic data from the application
of the Heisenberg flows associated with
$H_\n(y)$ and an element $y\in L$ of order $-1$.

In order to deal with Jacobian varieties, we further
assume that the field $k$ is the field $\mathbb{C}$ of
complex numbers in what follows in this section.
The computation (\ref{5.3}) shows that
every element of $H^1(C_\n,\mathcal{O}_{C_\n})$ 
is represented by
\begin{equation}
\sum_{j=1}^\ell\sum_{i=-\infty}^\infty t_{ij} y_j^{-i}\in
 \bigoplus_{j=1}^\ell \mathbb{C}((y_j)) = H_\n(y)\;.\label{5.4}
\end{equation}
The Heisenberg flows at $W$ are given by the equations
\begin{equation}
\frac{\partial W}{\partial t_{ij}} = y_j^{-i}\cdot W
= \big(h_{n_j}(y)\big)^{-i}\cdot W\;, \label{5.5}
\end{equation}
where $h_{n_j}(y)$ acts on $W$ through the block matrix
$$
\begin{pmatrix}
0\\
&\ddots\\
&&h_{n_j}(y)\\
&&&\ddots\\
&&&&0
\end{pmatrix}\;,
$$
and the index $i$ runs over all of $\mathbb{Z}$.
The formal integration
\begin{equation}
W(t) = \exp\left(\sum_{j=1}^\ell\sum_{i=-\infty}^\infty t_{ij} 
y_j^{-i}\right)
\cdot W\label{5.6}
\end{equation}
of the system (\ref{5.5}) shows that the stabilizers $A_0$
and $A_\n$ of $W(t)$ do not deform as $t$ varies, because 
the exponential factor
\begin{equation}
e(t) = \exp\left(\sum_{j=1}^\ell\sum_{i=-\infty}^\infty t_{ij} 
y_j^{-i}\right)\label{5.7}
\end{equation}
commutes with the algebra $H_\n(y)$. Note that half of the
 exponential factor
$$
\exp\left(\sum_{j=1}^\ell\sum_{i=-\infty}^{-1} t_{ij} 
y_j^{-i}\right)
$$
is an element of $\Gam_\n(y)$.

\medskip
\begin{thm} \label{5.8. Theorem}
 Let $y\in L$ be  a monic element of
order $-1$ and 
  $X_\n$  a finite type orbit of the Heisenberg
flows on $Gr_n(\mu)$ associated with 
$H_\n(y)$ starting at $W$.
 As we have seen in Proposition~\ref{5.1. Proposition}, the orbit $X_\n$
gives rise to a set of geometric data
$$
\left\langle f:\big( C_\n, \Delta, \Pi, \mathcal{F}, \Phi\big)
\longrarrow \big({C_0}, p, \pi, f_*\mathcal{F}, \phi\big)\right\rangle\;.
$$
Then the projection image $\overline{X}_\n$ of this orbit 
by $Q_{\n,y}:Gr_\n(\mu)\longrarrow Z_\n(\mu,y)$ is canonically
isomorphic to the Jacobian 
variety $\Jac(C_\n)$ of the curve $C_\n$ with $\overline{W}
= Q_{\n,y}(W)$ as its origin. Moreover, the orbit $\overline{X}_0$
of the KP system (written in terms of the variable $y$) defined on the
Grassmannian quotient $Z_\n(\mu,y)$ is isomorphic to the deformation
space
$$
\big\{ \mathcal{N}\tensor f_*\mathcal{F}\;\big|\; \mathcal{N}\in \Jac(C_0)\big\}\;.
$$
Thus we have a finite covering $\Jac(C_0)\longrarrow 
\overline{X}_0$ of the orbit, 
which is indeed isomorphic if $f_*\mathcal{F}$ is a
general vector bundle on $C_0$.
\end{thm}
\medskip

\begin{proof} 
Even though the formal integration (\ref{5.6}) is not well-defined as 
a point of the Grassmannian, we can still apply the same
construction of Section~\ref{sec: inverse} to the algebraic data 
$\big(A_0,A_\n,{W}(t)\big)$ understanding that the exponential 
matrix $e(t)$ of (\ref{5.7}) is an extra factor of degree 0. Of course
the curves, points, 
and the covering morphism $f:C_\n\longrarrow C_0$ are
the same as before. Therefore, we obtain
$$
\left\langle f:\big( C_\n, \Delta, \Pi, \mathcal{F}(t), \Phi(t)\big)
\longrarrow 
\big({C_0}, p, \pi, f_*\mathcal{F}(t), \phi(t)\big)\right\rangle\;,
$$
where the line bundle $\mathcal{F}(t)$
comes from the $A_\n$-module $W(t)$. 
We do not need to specify the data $\Phi(t)$ and $\phi(t)$
here, because they will disappear anyway by the trick of
Corollary~\ref{3.19. Corollary}. On the curve $C_\n$,
the formal expression $e(t)$ makes sense because of the homomorphism
$$
\exp: H^1(C_\n, \mathcal{O}_{C_\n})\owns 
\sum_{j=1}^\ell\sum_{i=-\infty}^\infty t_{ij} 
y_j^{-i}\longmapsto \big[e(t)\big]=\mathcal{L}(t)\in 
\Jac(C_\n)\subset H^1(C_\n,\mathcal{O}^*_{C_\n})\;,
$$
where $\mathcal{L}(t)$ is the line bundle of degree 0 corresponding to
the cohomology class
$\big[e(t)\big]\in  H^1(C_\n,\mathcal{O}^*_{C_\n})$. Thus the sheaf we obtain
from $W(t) = e(t)\cdot W$ is $\mathcal{F}(t) = \mathcal{L}(t)\tensor \mathcal{F}$.
Now consider the projection image $\big(A_0,A_\n,\overline{W}(t)\big)$
of the algebraic data by $Q_{\n,y}$. Then it corresponds to the data
\begin{equation}
\left\langle f:\bigg( C_\n, \Delta, \Pi, \mathcal{L}(t)\tensor
\mathcal{F}\bigg)
\longrarrow \bigg({C_0}, p, \pi, f_*\big(\mathcal{L}(t)\tensor
\mathcal{F}\big)\bigg)\right\rangle
\label{5.9}
\end{equation}
by Corollary~\ref{3.19. Corollary}. Since $\exp:H^1(C_\n, \mathcal{O}_{C_\n})
\longrarrow \Jac(C_\n)$ is surjective, we can 
define a map assigning (\ref{5.9}) to every point
$\mathcal{L}(t)\in \Jac(C_\n)$ of the Jacobian. 
Through the cohomology
functor, it gives indeed the desired identification of $\Jac(C_\n)$
and the orbit $\overline{X}_\n$:
$$
\Jac(C_\n)\owns \mathcal{L}(t)\longmapsto (\ref{5.9}) \longmapsto \overline{W}(t)
\in \overline{X}_\n\;.
$$
The KP system in the $y$-variable at $\overline{W}\in Z_\n(\mu,y)$
 is given by the equation
$$
\frac{\partial \overline{W}}{\partial s_m} = y^{-m}\cdot \overline{W}\;.
$$
The formal integration
$$
\overline{W}(s) = \exp\left( \sum_{m=1}^\infty s_my^{-m}\right)
\cdot \overline{W}
$$
corresponds to 
$$
\left\langle f:\bigg( C_\n, \Delta, \Pi, \big(f^*\mathcal{N}(s)\big)\tensor
\mathcal{F}\bigg)
\longrarrow \bigg({C_0}, p, \pi, \mathcal{N}(s)\tensor
f_*\mathcal{F}\bigg)\right\rangle\;,
$$
where $\mathcal{N}(s)\tensor f_*\mathcal{F}$ is the vector bundle 
corresponding to the $A_0$-module $\overline{W}(s)$. {From} (\ref{5.2}),
we have a surjective map of $H^1(C_0,\mathcal{O}_{C_0})$
onto the Jacobian variety $\Jac(C_0)
\subset H^1(C_0,\mathcal{O}^*_{C_0})$ defined by
$$
\exp: H^1(C_0,\mathcal{O}_{C_0})\owns \sum_{m=1}^\infty s_my^{-m}
\longmapsto \left[\exp\bigg(\sum_{m=1}^\infty s_my^{-m}\bigg)
\right] = \mathcal{N}(s)\in \Jac(C_0)\;.
$$
Thus the orbit $\overline{X}_0$ coincides with the deformation
space $\mathcal{N}(s)\tensor f_*\mathcal{F}$, which is covered
by $\Jac(C_0)$. The last statement of the theorem follows from
a result of \cite{L}.
This completes the proof.
\end{proof}

\medskip\noindent 
Let $(\eta_1, \cdots, \eta_\ell)$ be the transition function of 
$\mathcal{F}$ defined
on $U_j\setminus \{p_j\}$, where $\eta_j\in \mathbb{C}((y_j))$.
Then the family $\mathcal{F}(t)$ of line bundles on $C_\n$ is given by 
the transition function
$$
\left(\exp\bigg(\sum_{i=1}^\infty t_{i1} y_j^{-i}\bigg)\cdot\eta_1,\;
\cdots,\; 
\exp\bigg(\sum_{i=1}^\infty t_{i\ell} y_\ell^{-i}\bigg)\cdot\eta_\ell
\right)\;,
$$
and similarly, the line bundle $\mathcal{L}(t)$ is given by
$$
\left(\exp\bigg(\sum_{i=1}^\infty t_{i1} y_j^{-i}\bigg),\;
\cdots,\; 
\exp\bigg(\sum_{i=1}^\infty t_{i\ell} y_\ell^{-i}\bigg)
\right)\;.
$$
Here, we note that the nonnegative powers of $y_j=h_{n_j}(y)$ do not
contribute to these transition functions.

Recall that $H_\n(y)_0$
denotes the subalgebra of $H_\n(y)$ consisting of the
traceless elements. 

\medskip
\begin{thm} \label{5.10. Theorem}
In the same situation as above, the projection image 
$\overline{X}\subset Z_\n(\mu,y)$ of the orbit $X$ of 
the traceless Heisenberg flows ${\Psi}\big(H_\n(y)_0\big)$
starting at $\overline{W}$ is canonically isomorphic to
the Prym variety associated with the covering morphism 
$f:C_\n\longrarrow {C_0}$. 
\end{thm}
\medskip

\begin{proof}
Because of Remark~\ref{{1.3.} Remark}, the locus of $\mathcal{L}(t)\in \Jac(C_\n)$
such that 
$$
\det\bigg(f_*\big(\mathcal{L}(t)
\tensor \mathcal{F}\big)\bigg)= \det(f_*\mathcal{F})
$$
is the Prym variety $\Prym(f)$ associated with
 the covering morphism $f$.
So let us compute the factor
\begin{equation}
\mathcal{D}(t) = \det\big(f_*(\mathcal{L}(t)\otimes \mathcal{F})\big)
\otimes \det(f_*\mathcal{F})^{-1}\;,\label{5.11}
\end{equation}
which is a line bundle of degree 0 defined on $C_0$.
We use
the transition function $\eta$ of $f_*\mathcal{F}$ defined 
on $U_p\setminus \{p\}$ written in terms of the basis (\ref{3.6}). 
Since $f_*\mathcal{F}(t)$ is defined by
the $A_0$-module structure of $W(t) = e(t)\cdot W$, its transition
function is given by
$$
\exp
\begin{pmatrix}
\sum_{i=1}^\infty t_{i1} \big(h_{n_1}(y)\big)^{-i}\\
&\sum_{i=1}^\infty t_{i2} \big(h_{n_2}(y)\big)^{-i}\\
&&\ddots\\
&&&\sum_{i=1}^\infty t_{i\ell} \big(h_{n_\ell}(y)\big)^{-i}
\end{pmatrix}
\cdot \eta \;,
$$
where the $n\times n$ matrix  acts on the $y^{\alpha/n_j}$-part
of the basis of (\ref{3.6}) in an obvious way. 
Let us denote the above matrix by
$$
T(t) = 
\begin{pmatrix}
\sum_{i=1}^\infty t_{i1} \big(h_{n_1}(y)\big)^{-i}\\
&\sum_{i=1}^\infty t_{i2} \big(h_{n_2}(y)\big)^{-i}\\
&&\ddots\\
&&&\sum_{i=1}^\infty t_{i\ell} \big(h_{n_\ell}(y)\big)^{-i}
\end{pmatrix}\;.
$$
Then, it is clear that 
$\mathcal{D}(t) \isom \big[\exp \;\trace T(t)\big] 
\in H^1(C_0,\mathcal{O}^*_{C_0})$. 
{From} this expression, we see that if $\mathcal{L}(t)$ stays
on the orbit $\overline{X}$ of the traceless Heisenberg flows, then
 $\mathcal{D}(t) \isom \mathcal{O}_{C_0}$. Namely, $\overline{X}\subset \Prym(f)$.

Conversely, take a point $\overline{W}(t)\in \overline{X}_\n$
of the orbit of the Heisenberg flows defined on the quotient 
Grassmannian $Z_n(\mu,y)$. It corresponds to a unique
element $\mathcal{L}(t)\in \Jac(C_\n)$ by Theorem~\ref{5.8. Theorem}. Now suppose that 
the factor $\mathcal{D}(t)$
of (\ref{5.11}) is the trivial
bundle on $C_0$. Then it implies that $\big[\trace T(t)\big] = 0$
as an element of $H^1(C_0, \mathcal{O}_{C_0})$. In particular, 
$\trace T(t)$
acts on $\overline{W}$  trivially from (\ref{5.2}).  Therefore, 
$\overline{W}(t)$ is on the orbit of the flows defined by
$$
T(t) - I_n\cdot \frac{1}{n}\; \trace T(t)\;,
$$
which are clearly traceless. In other words, $\overline{W}(t)
\in \overline{X}$. Thus $\Prym(f)\subset \overline{X}$.
This completes the proof.
\end{proof}
\medskip
\begin{rem} \label{{5.12.} Remark}
Let us observe the case when  the curve $C_0$ downstairs
happens to be a $\mathbb{P}^1$. First of all, we note that the
$r$-reduced KP system associated with $y$
is nothing but the trace part of the Heisenberg flows 
defined by $H_\n(y)$. 
Because of the second half statement of Theorem~\ref{5.8. Theorem}, 
the trace part of the Heisenberg flows acts on the point
$\overline{W}\in Z_\n(\mu,y)$ trivially. Therefore, the 
 orbit $\overline{X}_\n$ of the entire Heisenberg flows
coincides with the orbit $\overline{X}$ of the traceless
part of the flows. Of course, this reflects the fact
that every Jacobian variety is a Prym variety associated with 
a covering over $\mathbb{P}^1$. Thus the
characterization theorem of Prym varieties we are
presenting below contains the characterization of Jacobians
of \cite{M1} as a special case.
\end{rem}

\medskip
Now consider the most trivial maximal commutative algebra 
$H = H_{(1, \cdots, 1)}(z) = \mathbb{C}((z))^{\dsum n}$. We define
the group $\Gam_{(1, \cdots, 1)}(z)$ following (\ref{2.15}), and denote by
\begin{equation}
Z_n(\mu) = Z_{(1, \cdots, 1)}(\mu,z) = 
Gr_n(\mu)\big/\Gam_{(1, \cdots, 1)}(z)\label{5.13}
\end{equation}
the corresponding 
Grassmannian quotient. On this space the algebra $H$ acts, 
and gives the $n$-component KP system. Let $H_0$ be the
traceless subalgebra of $H$, and consider the traceless 
$n$-component KP system on the Grassmannian quotient $Z_n(\mu)$.

\medskip
\begin{thm} \label{5.14. Theorem}
Every finite-dimensional orbit of the traceless $n$-component
KP system defined on the Grassmannian quotient 
$Z_n(\mu)$ of $(\ref{5.13})$ is canonically isomorphic to a (generalized) 
Prym variety. Conversely, every Prym variety associated with
a degree $n$ covering morphism of  smooth curves
can be realized in this way.
\end{thm}
\medskip

\begin{proof} The first half part has been already proved. So
start with the Prym variety $\Prym(f)$ associated with
a degree $n$ covering morphism $f:C\longrarrow C_0$ of 
smooth curves. Without
loss of generality, we can assume that $C_0$ is connected.
 Choose a point $p$ of $C_0$ outside of the branching locus so that
its preimage $f^{-1}(p)$ consists of $n$ distinct points
of $C$, and supply the necessary geometric objects to make
the situation into the geometric data
$$
\left\langle f:\big( C_\n, \Delta, \Pi, \mathcal{F}, \Phi\big)
\longrarrow \big({C_0}, p, \pi, f_*\mathcal{F}, \phi\big)\right\rangle
$$
of Definition~\ref{3.1. Definition} of rank one and type $\n = (1,\cdots,1)$
with $C= C_\n$.
 The data give rise to a unique triple
$(A_0,A_\n,W)$ of algebraic data by the cohomology
functor. We can choose $\pi = id$ 
so that the maximal commutative subalgebra we have
here is indeed $H = H_{(1, \cdots, 1)}(z)$. 
Define
$A'_0 = \{a\in \mathbb{C}((z))\;|\; a\cdot W\subset W\}$
and $A'= \{h\in H\;|\; h\cdot W\subset W\}$,
which satisfy $A_0\subset A'_0$ and $A_\n\subset A'$, and both
have  finite codimensions in the larger algebras. {From} the
triple of the algebraic data $(A'_0,A',W)$, we obtain
a set of geometric data 
$$
\left\langle {f'}:\big( C'_\n, \Delta', \Pi', \mathcal{F}', \Phi'\big)
\longrarrow \big({C'_0}, p', \pi', {f'}_*\mathcal{F}', \phi'\big)
\right\rangle\;.
$$
The morphism $(\alp, \beta, id)$ between the two sets of data 
consists of a morphism
$\alpha:C'_0\longrarrow C_0$
of the base curves and $\beta:C'\longrarrow C_\n$.
Obviously, these morphisms are birational, and hence, they have
to be an isomorphism, because $C_0$ and $C_\n$ are 
smooth. Going back to the algebraic data by the cohomology
functor, we obtain
$A_0 = A'_0$ and $A_\n = A'$. Thus the orbit of the 
traceless $n$-component KP system starting at $\overline{W}$
defined on the Grassmannian quotient $Z_n(\mu)$ is indeed the Prym
variety of the covering morphism $f$. This completes the proof of the
characterization theorem.
\end{proof}
\medskip
\begin{rem} \label{{5.15.} Remark}
In the above proof, we need the full information of the functor, 
not just the set-theoretical bijection of the objects. We use 
a similar argument once again in Theorem~\ref{6.15. Theorem}.
\end{rem}
\medskip
\begin{rem} \label{{5.16.} Remark}
The determinant line bundle $DET$ over $Gr_n(0)$ is defined by
$$
DET_W = \left(\bigwedge^{\max} \Ker\big(\gam_W\big)\right)^*
\bigotimes\bigwedge^{\max} \Coker\big(\gam_W\big)\;.
$$
The canonical section of the $DET$ bundle defines the determinant
divisor $Y$ of $Gr_n(0)$, whose support is the complement of
the big-cell $Gr^+_n(0)$. Note that the action of $\Gam_\n(y)$
preserves the big-cell. So we can define the \emph{big-cell}
of the Grassmannian quotient by $Z^+_\n(0, y) = 
Gr^+_n(0)\big/\Gam_\n(y)$. The determinant divisor also descends
to a divisor $Y/\Gam_\n(y)$, which we also call the determinant divisor
of the Grassmannian quotient. Consider a point $W\in Gr_n(0)$ at which
the Heisenberg flows of rank one produce a finite type 
orbit $X_\n$. The geometric data corresponding to this situation
consists of a curve $C_\n$ of genus $g = \dim_\mathbb{C} \overline{X}_\n$ 
and a line 
bundle $\mathcal{F}$ of degree $g - 1$ because of the Riemann-Roch formula
$$
\dim_\mathbb{C} H^0(C_\n,\mathcal{F}) - \dim_\mathbb{C} H^1(C_\n,\mathcal{F}) 
= \deg(\mathcal{F}) - r(g-1)\;.
$$
Thus we have an \emph{equality} $\overline{X}_\n = \Pic^{g-1}(C_\n)$ 
from the proof of Theorem~\ref{5.8. Theorem}. The intersection of $\overline{X}_\n$ with
the determinant divisor of $Z_\n(0,y)$ coincides with the theta
divisor $\Theta$ which gives the principal 
polarization of $\Pic^{g-1}(C_\n)$. However, the restriction
of this divisor to the Prym variety  
does not give a principal polarization 
as we have noted in Section~\ref{sec: 1}.
\end{rem}
\medskip
\begin{rem} \label{{5.17.} Remark}
{From} the expression of (\ref{5.9}), we can see that
a finite-dimensional orbit of the
Heisenberg flows of rank one defined on the Grassmannian quotient
gives a family of deformations $f_*\big(\mathcal{L}(t)\tensor \mathcal{F}\big)$
of the vector bundle $f_*\mathcal{F}$ on $C_0$. It is an interesting 
question to ask what kind of deformations does this family produce. 
More generally, we can ask the following question:
For a given curve and a family of vector bundles on it, can one find
a point $W$ of the Grassmannian $Gr_n(\mu)$ and a suitable
Heisenberg flows such that the orbit starting from $W$
contains the original family?

It is known that for every vector bundle $\mathcal{V}$ of rank $n$
on a smooth curve $C_0$,
there is a degree $n$ covering $f:C\longrarrow C_0$ and a line
bundle $\mathcal{F}$ on $C$ such that $\mathcal{V}$ is
isomorphic to the direct image sheaf $f_*\mathcal{F}$. We can supply
 suitable local data so that we have a set of
geometric data
 $$
\left\langle f:\big(C_\n, \Delta, \Pi, \mathcal{F}\big)
\longrarrow \big({C_0}, p, \pi, f_*\mathcal{F}\big)\right\rangle
$$
with $C_\n = C$. Let $(A_0, A_\n, \overline{W})$ be the
triple of algebraic data  corresponding to the 
above geometric situation
with a point $\overline{W}\in Z_\n(\mu,z)$, 
where $\mu$ is the Euler characteristic of the original bundle $\mathcal{V}$.
 Now the problem
is to compare the family of deformations given by (\ref{5.9}) and the
original family. 

The only thing we can say about this question at the present moment
is the following. 
If the original vector bundle is a general stable bundle,
then one can find a set of
 geometric data  and a corresponding point $\overline{W}$
of a Grassmannian quotient such that 
 there is a dominant and generically finite map of
 a Zariski open subset of the orbit of the Heisenberg flows
starting from  $\overline{W}$ into
the moduli space of stable vector bundles of rank $n$ and
degree $\mu + n\big(g(C_0) - 1\big)$ over the curve $C_0$.
Note that this
statement is just an interpretation of a theorem
of \cite{BNR} into our language using Theorem~\ref{5.8. Theorem}. 

As in the proof of Theorem~\ref{5.14. Theorem}, the Heisenberg flows can be replaced
by the $n$-component KP flows if we choose the  point 
$p\in C_0$ away from
the branching locus of $f$. Thus one may say that the $n$-component
KP system can produce general vector bundles of rank $n$ 
defined on an arbitrary smooth curve in its orbit.
\end{rem}
\bigskip

\section{Commuting ordinary differential operators with matrix 
coefficients}
\label{sec: ODE}
\medskip

In this section, we work with an arbitrary field $k$ again.
Let us denote by 
$$
E = \big(k[[x]]\big)((\partial^{-1}))
$$
the set of all 
\hyphenation{pseu-do-dif-fer-en-tial}
pseudodifferential operators with coefficients in 
$k[[x]]$, where $\partial = d/dx$. This is an associative algebra
and has a natural filtration 
$$
E^{(m)} = \big(k[[x]]\big)[[\partial^{-1}]]\cdot \partial^m
$$ 
defined by the \emph{order} of the operators.
We can identify $k((z))$ 
with the set of pseudodifferential operators with 
constant coefficients by the
\emph{Fourier transform} $z$ $=$ $\partial^{-1}$:
$$
L = k((z)) = k((\partial^{-1}))\subset E\;.
$$
There is also a canonical projection
\begin{equation}
\rho: E \longrarrow E/Ex \isom k((\partial^{-1})) = L\;,\label{6.1}
\end{equation}
where $Ex$ is the left-maximal ideal of $E$ generated by $x$.
In an explicit form, this projection is given by 
\begin{equation}
\rho : E\owns P = \sum_{m\in \mathbb{Z}} \partial^m\cdot a_m(x)
\longmapsto \sum_{m\in\mathbb{Z}} a_m(0)z^{-m}\in L\;.\label{6.2}
\end{equation}
It is obvious from (\ref{6.1}) that $L$ is a left $E$-module. The action 
is given by $P\cdot v = P\cdot\rho(Q) = \rho(PQ)$, where
$v\in L = E/Ex$ and $Q\in E$ is a representative of the equivalence
class such that $\rho(Q) = v$. The well-definedness of this action
is easily checked.
We also use the notations
$$
\begin{cases}
D = \big(k[[x]]\big)[\partial]\\
E^{(-1)} = \big(k[[x]]\big)[[\partial^{-1}]]\cdot\partial^{-1}\;,
\end{cases}
$$
which are the set of linear ordinary differential operators
and the set of pseudodifferential operators of negative order,
respectively. Note that there is a natural 
left $\big(k[[x]]\big)$-module direct sum decomposition
\begin{equation}
E = D\dsum E^{(-1)}\;.\label{6.3}
\end{equation}
According to this decomposition, we write $P = P^+ \dsum P^-$, 
$P \in E$, $P^+ \in D$, and $P^-\in E^{(-1)}$.

Now consider the matrix algebra $gl(n,E)$ defined over the 
noncommutative algebra $E$, which is the algebra of pseudodifferential
operators with coefficients in matrix valued functions. 
This algebra acts on our vector space
$V = L^{\dsum n} \isom \big(E/Ex\big)^{\dsum n}$ from the left.
In particular, every element of $gl(n,E)$ gives rise to a
vector field on the Grassmannian $Gr_n(\mu)$ via (\ref{2.4}).
The decomposition (\ref{6.3}) induces 
$$
V = k[z^{-1}]^{\dsum n}\dsum \big(k[[z]]\cdot z\big)^{\dsum n}
$$
after the identification $z = \partial^{-1}$, 
and the base point $k[z^{-1}]^{\dsum n}$ of the Grassmannian
$Gr_n(0)$ of index 0 is the residue class of $D^{\dsum n}$
in $E^{\dsum n}$ via the projection $E^{\dsum n}
\longrarrow E^{\dsum n}\big/\big(E^{(-1)})^{\dsum n}$. Therefore,
the $gl(n,D)$-action on $V$ preserves $k[z^{-1}]^{\dsum n}$.
The following proposition shows that the converse is also true:

\medskip
\begin{prop} \label{6.4. Proposition}
A pseudodifferential operator $P\in gl(n,E)$ with matrix coefficients
is a differential operator, i.e.~$P\in gl(n,D)$, if and only if
$$
P\cdot k[z^{-1}]^{\dsum n}\subset k[z^{-1}]^{\dsum n}\;.
$$
\end{prop}
\medskip
\begin{proof}
The case of $n = 1$ of this proposition was established in
Lemma~7.2 of \cite{M3}. So let us assume that 
$P = \big(P_{\mu\nu}\big)\in gl(n,E)$ preserves the base point
$k[z^{-1}]^{\dsum n}$. If we apply the matrix $P$ to the vector
subspace 
$$
0\dsum\cdots\dsum 0\dsum k[z^{-1}]\dsum 0\dsum \cdots \dsum 0
\subset k[z^{-1}]^{\dsum n}
$$ 
with only nonzero entries in the $\nu$-th position, 
then we know that $P_{\mu\nu}\in E$ stabilizes $k[z^{-1}]$ in $L$. 
Thus $P_{\mu\nu}$ is a differential operator, i.e.\ $P\in gl(n,D)$.
This completes the proof.
\end{proof}
\medskip
Since  differential operators preserve the base point
of the Grassmannian $Gr_n(0)$, the negative order pseudodifferential
operators should give the most part of $Gr_n(0)$. In fact, we have
\medskip
\begin{thm} \label{6.5. Theorem}
Let $S\in gl(n,E)$ be a monic zero-th order pseudodifferential
operator of the form
\begin{equation}
S = I_n + \sum_{m = 1}^\infty s_m(x) \partial^{-m} \;,\label{6.6}
\end{equation}
where $s_m(x)\in gl\big(n,k[[x]]\big)$. Then the map
$$
\sigma : \Sigma\owns S\longmapsto W = S^{-1}\cdot k[z^{-1}]^{\dsum n}
\in Gr^+_n(0)
$$
gives a bijective correspondence between the set $\Sigma$ of
pseudodifferential operators of the form of $(\ref{6.6})$
and the big-cell $Gr_n^+(0)$ of the index 0 Grassmannian.
\end{thm}
\medskip
\begin{proof} Since $S$ is invertible of order 0, 
we have $S^{-1}\cdot V = V$
and $S^{-1}\cdot V^{(-1)} = V^{(-1)}$, where
$V^{(-1)} = F^{(-1)}(V) = \big(k[[z]]z\big)^{\dsum n}$.
Thus $V = S^{-1}\cdot
k[z^{-1}]^{\dsum n} \dsum V^{(-1)}$, which shows that 
$\sigma$ maps into the big-cell. 

The injectivity of $\sigma$ is easy: if $S_1^{-1}\cdot 
k[z^{-1}]^{\dsum n} = S_2^{-1}\cdot k[z^{-1}]^{\dsum n}$, then
$S_1S_2^{-1}\cdot k[z^{-1}]^{\dsum n} = k[z^{-1}]^{\dsum n}$. 
It means, by Proposition~\ref{6.4. Proposition}, that $S_1S_2^{-1}$ is a differential
operator. Since $S_1S_2^{-1}$ has the same form of (\ref{6.6}), the
only possibility is that $S_1S_2^{-1} = I_n$, which implies the
injectivity of $\sigma$.

In order to establish surjectivity, take an arbitrary point $W$ 
of the big-cell $Gr^+(0)$. We can choose a basis $\big\langle
\w_j^\mu\big\rangle_{1\le j\le n, 0\le \mu}$ for the vector
space $W$ in the form 
$$
\w_j^\mu = \e_jz^{-\mu} + \sum_{\nu = 1}^\infty \sum_{i=1}^n
\e_iw_{j\nu}^{i\mu}z^\nu\;,
$$
where $\e_j$ is the elementary column vector of size $n$
and $w_{j\nu}^{i\mu}\in k$. Our goal is to construct an
operator $S\in \Sigma$ such that $S^{-1}\cdot k[z^{-1}]^{\dsum n} = W$.
Let us put $S^{-1} = \big(S^i_j\big)_{1\le i,j\le n}$
with 
$$
S^i_j = \delta^i_j + \sum_{\nu=1}^\infty \partial^{-\nu} \cdot
s^i_{j\nu}(x)\;.
$$ 
Since every coefficient $s^i_{j\nu}(x)$
of $S^{-1}$ is a formal power series in $x$, we can construct
the operator by induction on the power of $x$. So let us
assume that we have constructed $s^i_{j\nu}(x)$ modulo $k[[x]]x^\mu$.
We have to introduce one more equation of order $\mu$ 
in order to determine the coefficient of $x^\mu$ in $s^i_{j\nu}(x)$,
which comes from the equation
$$
S^{-1}\cdot \e_jz^{-\mu} = {\text{ a linear combination of   }}
\w_i^{\nu}\;.
$$
For the purpose of finding
a consistent equation, 
let us compute the left-hand side by using the projection $\rho$
of (\ref{6.2}):

$$
\begin{aligned}
S^{-1}\cdot \e_jz^{-\mu} &= \sum_{i=1}^n \e_i\cdot S^i_j\cdot z^{-\mu}\\
&= \e_jz^{-\mu} + \rho\left(\sum_{\nu=1}^\infty\sum_{i=1}^n 
\partial^{-\nu}\cdot s^i_{j\nu}(x)\e_i\cdot \partial^\mu\right)\\
&= \e_jz^{-\mu} + \rho\left(\sum_{\nu=1}^\infty\sum_{i=1}^n
\sum_{m=0}^\mu (-1)^m\binom \mu{m} \e_i\cdot \partial^{\mu-\nu-m}
\cdot {s^i_{j\nu}}^{(m)}(x)\right)\\
&= \e_jz^{-\mu} + \sum_{\nu=1}^\infty\sum_{i=1}^n
\sum_{m=0}^\mu (-1)^m\binom \mu{m}  {s^i_{j\nu}}^{(m)}(0)\cdot
\e_i z^{-\mu+\nu+m}\\
&= \e_jz^{-\mu} + \sum_{\alp = 1}^{\mu-1}\sum_{m=0}^\alp \sum_{i=1}^n
(-1)^m\binom \mu{m}  {s^i_{j\alp - m}}^{(m)}(0)\cdot
 \e_i z^{-\mu+\alp}\\
&\qquad\qquad + \sum_{m=0}^{\mu-1} \sum_{i=1}^n
(-1)^m\binom \mu{m} {s^i_{j\mu - m}}^{(m)}(0)\cdot
  \e_i \\
&\qquad\qquad + \sum_{\beta=1}^\infty \sum_{m=0}^\mu\sum_{i=1}^n
(-1)^m \binom \mu{m} {s^i_{j\beta+\mu - m}}^{(m)}(0)
\cdot\e_i z^\beta\;.
\end{aligned}
$$
 
\noindent
Thus we see that the equation
\begin{equation}
\begin{aligned}
S^{-1}\cdot \e_jz^{-\mu} = \w_j^\mu
&+ \sum_{\alp = 1}^{\mu-1}\sum_{m=0}^\alp \sum_{i=1}^n
(-1)^m\binom \mu{m}  {s^i_{j\alp - m}}^{(m)}(0)\cdot
 \w_i^{\mu-\alp}\\
&+ \sum_{m=0}^{\mu-1} \sum_{i=1}^n
(-1)^m\binom \mu{m} {s^i_{j\mu - m}}^{(m)}(0)\cdot
  \w_i^0
\end{aligned} \label{6.7}
\end{equation}
is the identity for the coefficients of $\e_iz^{-\nu}$ for all 
$i$ and $\nu\ge 0$, and determines $s^i_{j\beta}(0)^{(\mu)}$
uniquely, because the coefficient of $s^i_{j\beta}(0)^{(\mu)}$
in the equation is $(-1)^\mu$. Thus by solving (\ref{6.7}) for 
all $j$ and $\mu\ge 0$ inductively, we can determine the operator
$S$ uniquely, which satisfies the desired property by the
construction.
This completes the
proof.
\end{proof}
\medskip
Using this identification of $Gr^+(0)$ and $\Sigma$, 
we can translate the Heisenberg flows defined on the big-cell
into a system of nonlinear partial differential equations.
Since we are not introducing any analytic structures in $\Sigma$,
we cannot talk about a Lie group structure in it. However, the 
exponential map
$$
\exp: gl\big(n,E^{(-1)}\big)\longrarrow I_n + gl\big(n,E^{(-1)}\big)
 = \Sigma
$$
is well-defined and surjective, and hence we can regard 
$gl\big(n,E^{(-1)}\big)$
as the \emph{Lie algebra} of the infinite-dimensional group $\Sigma$.
Symbolically, we have an identification
\begin{equation}
T_{k[z^{-1}]^{\dsum n}}Gr^+(0) \isom
gl\big(n,E^{(-1)}\big) = {\text{Lie}}(\Sigma) 
= T_{I_n}\Sigma = S^{-1}\cdot 
T_S\Sigma \label{6.8}
\end{equation}
for every $S\in \Sigma$.
The equation 
$$
\frac{\partial W(t)}{\partial t_{ij}} = 
\big(h_{n_j}(y)\big)^{-i}\cdot W(t) 
$$
is an equation of tangent vectors at the point $W(t)$. 
We now identify the variable $y$ of (\ref{2.7}) with a pseudodifferential
operator 
\begin{equation}
y = \partial^{-r} + \sum_{m=1}^\infty c_m \partial^{-r-m} \label{6.9}
\end{equation}
with coefficients in $k$.
Then the block matrix $h_{n_j}(y)$ 
of (\ref{5.5}) is identified with an element
of $gl(n,E)$. Let $W(t)$ be a solution of (\ref{5.5}) which lies in
$Gr^+_n(0)$, where $t = (t_{ij})$. Writing $W(t) = S(t)^{-1}\cdot
k[z^{-1}]^{\dsum n}$, the tangent vector of the left-hand 
side of (\ref{5.5}) is given by
$$
\frac{\partial W(t)}{\partial t_{ij}} = \frac{\partial
S(t)^{-1}}{\partial t_{ij}}\;,
$$
which then gives an element 
$$
S(t)\cdot \frac{\partial S(t)^{-1}}{\partial t_{ij}}
= -\; \frac{\partial S(t)}{\partial t_{ij}} \cdot S(t)^{-1}
\in E^{(-1)}
$$
by (\ref{6.8}).
The tangent vector of the right-hand side of (\ref{5.5}) 
is $\big(h_{n_j}(y)\big)^{-i}\in \Hom_{\text{cont}}(W,V/W)$, 
which gives rise to
a tangent vector $S(t)\cdot \big(h_{n_j}(y)\big)^{-i}\cdot S(t)^{-1}$
at the base point $k[z^{-1}]^{\dsum n}$ of the big-cell
by the diagram
$$
\begin{CD}
k[z^{-1}]^{\dsum n} @>>> V @>{S\cdot h\cdot S^{-1}}>> V 
@>>> V\big/k[z^{-1}]^{\dsum n}\\
@V{S^{-1}}VV @V{S^{-1}}VV @VV{S^{-1}}V @VV{S^{-1}}V\\
W @>>> V @>>{h}> V @>>> V/W,
\end{CD}
$$
where we denote $W = W(t)$, $S = S(t)$ and $h = 
\big(h_{n_j}(y)\big)^{-i}$.
Since the base point is preserved by the differential operators,
the equation of the tangent vectors reduces to an equation 
\begin{equation}
\frac{\partial S(t)}{\partial t_{ij}}\cdot S(t)^{-1}
= -\bigg(S(t)\cdot \big(h_{n_j}(y)\big)^{-i}\cdot S(t)^{-1}\bigg)^-
\label{6.10}
\end{equation}
in the Lie algebra $gl\big(n,E^{(-1)}\big)$ level,
where $(\bullet)^-$ denotes the negative order part of the 
operator by (\ref{6.3}).
We call this equation the \emph{Heisenberg KP system}. 
Note that the above equation is trivial for negative $i$
because of (\ref{6.9}). In terms of the operator
$$
P(t) = S(t)\cdot y^{-1}\cdot I_n\cdot S(t)^{-1}\in gl(n,E)
$$
whose leading term is $I_n\cdot \partial^r$, the equation (\ref{6.10})
becomes a more familiar \emph{Lax equation}
$$
\frac{\partial P(t)}{\partial t_{ij}} = 
\left[\bigg(S(t)\cdot \big(h_{n_j}(y)\big)^{-i}\cdot S(t)^{-1}\bigg)^+,
\;P(t)\right]\;.
$$
In particular, the Heisenberg KP system describes 
 infinitesimal isospectral
deformations of the operator $P = P(0)$.
Note that if one chooses $y=z=\partial^{-1}$ in (\ref{6.9}), then
the above Lax equation for the case of $n=1$ becomes the original KP 
system.
We can solve the initial value problem of the Heisenberg KP system
(\ref{6.10}) by the \emph{generalized Birkhoff decomposition} of \cite{M2}:
\begin{equation}
\exp\left(\sum_{j=1}^\ell\sum_{i=1}^\infty t_{ij} 
\big(h_{n_j}(y)\big)^{-i}\right)\cdot S(0)^{-1} = S(t)^{-1}\cdot
Y(t)\;,\label{6.11}
\end{equation}
where $Y(t)$ is an invertible differential operator of infinite
order defined in \cite{M2}. In order to see that the $S(t)$ of (\ref{6.11})
gives a solution of (\ref{6.10}), we differentiate the equation (\ref{6.11})
with respect to $t_{ij}$. Then we have
$$
S(t)\cdot \big(h_{n_j}(y)\big)^{-i}\cdot S(t)^{-1}
= -\;\frac{\partial S(t)}{\partial t_{ij}}\cdot S(t)^{-1}
+ \frac{\partial Y(t)}{\partial t_{ij}}\cdot Y(t)^{-1}\;,
$$
whose negative order terms are nothing but the 
Heisenberg KP system (\ref{6.10}). It shows that the Heisenberg KP system
is a completely integrable system of nonlinear partial differential
equations.

Now, consider a set of geometric data
$$
\left\langle f:\big( C_\n, \Delta, \Pi, \mathcal{F}, \Phi\big)\longrarrow \big({C_0}, p, \pi, f_*\mathcal{F}, \phi\big)\right\rangle
$$
such that $H^0(C_\n,\mathcal{F}) = H^1(C_\n,\mathcal{F}) = 0$.
Then by the cohomology functor of Theorem~\ref{3.17. Theorem}, it gives rise to a
triple $(A_0, A_\n, W)$ satisfying that $W\in Gr^+(0)$. By Theorem~\ref{6.5. Theorem},
there is a monic zero-th order pseudodifferential operator $S$
such that $W = S^{-1}\cdot k[z^{-1}]^{\dsum n}$. 
Using the identification
(\ref{6.9}) of the variable $y$ as the pseudodifferential operator
with constant coefficients, we can define two commutative subalgebras
of $gl(n,E)$ by
\begin{equation}
\begin{cases}
B_0 = S\cdot A_0\cdot S^{-1} \\
B_\n = S\cdot A_\n\cdot S^{-1}\;.
\end{cases}\label{6.12}
\end{equation}
The inclusion relation $A_0\subset k((y))$ gives us
$B_0\subset k((P^{-1}))$, where 
$P =$ $S\cdot y^{-1}\cdot I_n\cdot S^{-1}$ $\in$ $gl(n,E)$.
Since $A_0$ and $A_\n$ stabilize $W$, we know that $B_0$ and
$B_\n$ stabilize $k[z^{-1}]^{\dsum n}$. Therefore, these
algebras are commutative algebras of ordinary differential
operators with matrix coefficients!
\medskip
\begin{Def} \label{6.13. Definition}
We denote by $\mathcal{C}^+(\n,0,r)$ the set of objects
$$
\left\langle f:\big( C_\n, \Delta, \Pi, \mathcal{F}, \Phi\big)\longrarrow 
\big({C_0}, p, \pi, f_*\mathcal{F}, \phi\big)\right\rangle
$$
of the category $\mathcal{C}(\n)$ of index 0 and rank $r$ such that
$$
H^0(C_\n, \mathcal{F}) = H^1(C_\n, \mathcal{F}) = 0\;.
$$
The set of pairs $(B_0,B)$ of commutative algebras satisfying
the following conditions is denoted by $\mathcal{D}(n,r)$:
\begin{enumerate}
\item $k\subset B_0\subset B\subset gl(n,D)$.

\item $B_0$ and $B$ are commutative $k$-algebras.

\item There is an operator $P\in gl(n,E)$ whose leading term
is $I_n\cdot \partial^r$ such that $B_0\subset k((P^{-1}))$.

\item The projection map $B_0\longrarrow k((P^{-1}))\big/k[[P^{-1}]]$
is Fredholm.

\item $B$ has rank $n$ as a torsion-free module over $B_0$.
\end{enumerate}
\end{Def}
\medskip\noindent
Using this definition, we can summarize
\medskip
\begin{prop} \label{6.14. Proposition}
The construction $(\ref{6.12})$ gives a canonical map 
$$
\chi_{\n,r}: \mathcal{C}^+(\n,0,r) \longrarrow \mathcal{D}(n,r)
$$
for every $r$ and a positive 
integral vector $\n = (n_1,\cdots, n_\ell)$
with $n = n_1+\cdots  + n_\ell$.
\end{prop}
\medskip
\noindent
If the field $k$ is of characteristic zero, then we can construct
maximal commutative algebras of ordinary differential operators
with coefficients in matrix valued functions as an application of the
above proposition.
\medskip
\begin{thm} \label{6.15. Theorem}
Every set 
$$
\left\langle f:\big( C_\n, \Delta, \Pi, \mathcal{F}, \Phi\big)\longrarrow 
\big({C_0}, p, id, f_*\mathcal{F}, \phi\big)\right\rangle
$$
of geometric data with a smooth curve $C_\n$, $\pi = id$ 
and a line bundle
$\mathcal{F}$ satisfying that
 $H^0(C_\n,\mathcal{F}) = H^1(C_\n,\mathcal{F}) = 0$ 
gives rise to a maximal commutative subalgebra 
$B_\n\subset gl(n,D)$ by $\chi_{\n,1}$.
\end{thm}
\medskip
\begin{proof}
Let $(B_0,B_\n)$ be the image of $\chi_{\n,1}$ applied to the 
above object, and $(A_0,A_\n,W)$ the stabilizer data corresponding
to the geometric data. Recall that $B_0 = S\cdot A_0\cdot
S^{-1}$, where $S$ is the operator corresponding to $W$.
Since $r = 1$ in our case, (\ref{4.8})
implies the existence of
an element $a\in A_0$ of the form
$$
a = a(z^{-1}) = z^{-m} + c_2 z^{-m+2} + c_3 z^{-m+3} + \cdots
\in A_0\subset k((z))\;.
$$
We call a pseudodifferential operator $a(\partial)\cdot I_n
\in gl(n,E)$ a \emph{normalized} scalar
diagonal operator of order $m$ with constant coefficients.
Here, we need
\end{proof}
\medskip
\begin{lem} \label{6.16. Lemma}
Let $K\in gl(n,E)$ be a normalized scalar
diagonal operator of order $m>0$ with constant coefficients
and $Q = (Q_{ij})$ an arbitrary element of $gl(n,E)$. If $Q$ and $K$
commute, then every coefficient of $Q$ is a constant matrix.
\end{lem}
\medskip
\begin{proof} Let $K = a(\partial)\cdot I_n$ for some
$a(\partial)\in k((\partial^{-1}))$. It is well known that
there is a monic zero-th order pseudodifferential operator 
$S_0\in E$ such that 
$$
S_0^{-1}\cdot a(\partial)\cdot S_0 = 
\partial^m\;.
$$
Since $a(\partial)$ is a constant coefficient operator, we can show
that (see \cite{M3})
$$
S_0^{-1}\cdot k((\partial^{-1})) \cdot S_0 = k((\partial^{-1}))\;.
$$
Going back to the matrix case, we have
$$
0 = (S_0\cdot I_n)^{-1}\cdot [Q,K]\cdot (S_0\cdot I_n)
= \big[(S_0\cdot I_n)^{-1}\cdot Q \cdot (S_0\cdot I_n),\;\partial^m
\cdot I_n\big]\;.
$$
In characteristic zero, commutativity with $\partial^m$
implies commutativity with $\partial$. Thus each matrix
component $S_0^{-1}\cdot Q_{ij}\cdot S_0$ commutes with 
$\partial$, and hence 
$S_0^{-1}\cdot Q_{ij}\cdot S_0$ $\in$ $k((\partial^{-1}))$. 
Therefore, $Q_{ij}\in k((\partial^{-1}))$. This completes
the proof of lemma.
\end{proof}
\medskip\noindent
Now, let $B \supset B_\n$ be a commutative subalgebra
of $gl(n,D)$ containing $B_\n$.
Since $B_0 = S\cdot A_0\cdot S^{-1}$ and $B_0\subset B$,
every element of $B$ commutes with $S\cdot a(\partial)\cdot I_n\cdot
S^{-1}$. Then by the lemma, we have
$$
A = S^{-1}\cdot B\cdot S\subset gl\big(n,k((\partial^{-1}))\big)\;.
$$
Note that the algebra $A$ stabilizes $W = S^{-1}\cdot 
k[z^{-1}]^{\dsum n}$. Since $H_\n(z)$
can be  generated by $A_\n$ over $k((z)) = k((\partial^{-1}))$,
every element of $A$ commutes with $H_\n(z)$. Therefore,
 we have $A\subset H_\n(z)$ because of the maximality
of $H_\n(z)$.
Thus we obtain another triple $(A_0,A,W)$ of stabilizer data
of the same type $\n$.
The inclusion $A_\n\longrarrow A$ gives rise to a 
birational morphism
$\beta:C\longrarrow C_\n$. Since we are assuming that the curve $C_\n$ 
is nonsingular, $\beta$ has to be an isomorphism, which then
implies that $A = A_\n$. Therefore, we have $B = B_\n$.
This completes the proof of maximality of $B_\n$.

\medskip

\begin{rem} \label{{6.17.} Remark}
There are other maximal commutative subalgebras in $gl(n,D)$ than
what we have constructed in Theorem~\ref{6.15. Theorem}. It corresponds to the 
fact that the algebras $H_\n(z)$ are not the only maximal commutative
subalgebras of the formal loop algebra $gl\big(n, k((z))\big)$.
\end{rem}

\bigskip

\providecommand{\bysame}{\leavevmode\hbox to3em{\hrulefill}\thinspace}

\end{document}